\documentclass[11pt,letterpaper]{amsart}
\usepackage[utf8]{inputenc}
\usepackage[T1]{fontenc}
\input{glyphtounicode}
\usepackage[english]{babel}
\usepackage{amsmath,amssymb,amsthm}
\usepackage{array}
\usepackage{xcolor}

\usepackage[accsupp]{axessibility}
\usepackage{hyperref}
\hypersetup{
  hidelinks,
  pdflang={en-US},
  pdftitle={Only constant mean curvature surfaces in E(kappa,tau), tau not 0, have a holomorphic original Abresch-Rosenberg differential},
  pdfauthor={Hilario Alencar and Harold Rosenberg},
  pdfsubject={Differential geometry: surfaces in homogeneous 3-manifolds},
  pdfkeywords={Abresch-Rosenberg differential; Hopf differential; constant mean curvature; homogeneous 3-manifolds; Puiseux series},
  pdfdisplaydoctitle=true,
  bookmarksnumbered=true
}

\theoremstyle{plain}
\newtheorem{theorem}{Theorem}[section]
\newtheorem{proposition}[theorem]{Proposition}
\newtheorem{lemma}[theorem]{Lemma}
\newtheorem{corollary}[theorem]{Corollary}
\newtheorem{fact}[theorem]{Fact}
\theoremstyle{definition}
\newtheorem{remark}[theorem]{Remark}

\numberwithin{equation}{section}
\newcommand{\R}{\mathbb R}
\newcommand{\C}{\mathbb C}
\newcommand{\Q}{\mathbb Q}
\newcommand{\Z}{\mathbb Z}

\newcommand{\E}{E(\kappa,\tau)}
\newcommand{\QAR}{Q_{\mathrm{AR}}}
\newcommand{\calQAR}{\mathcal Q_{\mathrm{AR}}}

\title[Holomorphic original AR differential and CMC]{Only constant mean curvature surfaces in $E(\kappa,\tau)$, $\tau\ne0$, have a holomorphic original Abresch--Rosenberg differential}
\author{Hil\'ario Alencar}
\address{Instituto de Matem{\'a}tica, Universidade Federal de Alagoas, 57072-900, Macei\'o, AL, Brazil}
\email{hilario@mat.ufal.br}
\author{Harold Rosenberg}
\address{13 rue du Faubourg Montmartre, 75009 Paris, France}
\email{rosenbergharold@gmail.com}
\thanks{\emph{Use of generative AI.} Claude (Anthropic) was used for English-language editing, for the typesetting of the tables, whose content was produced by the authors, and for adapting the \LaTeX{} source to the journal's template. It was also used to discuss and check selected calculations. All mathematical content, statements, proofs, and references are the authors' own, and the authors take full responsibility for them.}
\thanks{H.~Alencar was partially supported by the Brazilian National Council for Scientific and Technological Development -- CNPq, Brazil (grant 303118/2022-9), and by the Coordination for the Improvement of Higher Education Personnel -- CAPES, Brazil (Finance Code 001).}

\begin{document}
\begin{abstract}
For immersed surfaces in the homogeneous spaces $E(\kappa,\tau)$ with
$\tau\ne0$, we prove that the original Abresch--Rosenberg differential
$[2(H+i\tau)p-(\kappa-4\tau^2)A^2]\,d\zeta^2$ is holomorphic if and only if
the mean curvature $H$ is constant; here $p\,d\zeta^2$ is the Hopf
differential and $A=\langle\xi,\partial_\zeta\rangle$. The result is local:
no completeness, compactness or analyticity is assumed, and there are no
exceptional surfaces. Where $dH\ne0$, holomorphy forces the shape operator into the normal form
$\bigl(\begin{smallmatrix}0&\tau\\\tau&2H\end{smallmatrix}\bigr)$ in the frame
$(\nabla H/|\nabla H|,J\nabla H/|\nabla H|)$. The round case $\kappa=4\tau^2$ is
settled by an elementary argument; for $\kappa\ne4\tau^2$ the structure
equations become a rational system whose solutions are algebraic, and an
embedding into a field of Puiseux series excludes every branch. In $M^2(\kappa)\times\R$,
$\kappa\ne0$, where the statement fails, we show that a non-CMC surface with
holomorphic original differential has nowhere vanishing $dH$ and constant
angle function, and belongs to an explicit family; in particular, closed
surfaces with holomorphic original differential are CMC there. We also show that
the Hopf-cylinder alternative in a theorem of Dorfmeister, Inoguchi and
Kobayashi contains no non-CMC surface.
\end{abstract}


\maketitle

\section{Introduction and main results}

Hopf's differential is holomorphic precisely on constant mean curvature
(CMC) surfaces in three-dimensional space forms; its vanishing on a
topological sphere gives Hopf's sphere theorem \cite[Part~Two, Ch.~VI, \S1.4, Lemma~2.2 and Theorems~2.1 and~2.6, pp.~137--140]{Hopf1989}.
Abresch and Rosenberg extended this construction to CMC surfaces in
$S^2\times\R$ and $\mathrm H^2\times\R$ \cite[(1)--(2) and Theorem~1, p.~142]{AR2004} and, more generally,
in $\E$ \cite[Theorems~4 and~5, pp.~16--17]{AR2005}.

For $\tau\ne0$, $\E$ is the simply connected homogeneous Riemannian
$3$-manifold fibering over $M^2(\kappa)$ by a Riemannian submersion with
geodesic fibers, whose unit vertical Killing field $\xi$ satisfies
$\overline\nabla_X\xi=\tau X\times\xi$.
The signs $\kappa<0$, $\kappa=0$, and $\kappa>0$ give
$\widetilde{PSL}_2(\R)$, $\mathrm{Nil}_3$, and Berger spheres,
respectively \cite[\S2, p.~91]{Daniel2007}; see also \cite{Scott1983}.
For $\tau=0$, $E(\kappa,0)$ is the product $M^2(\kappa)\times\R$ and
$\xi=\partial_t$ is parallel. Write $k=\kappa-4\tau^2$; for $\tau\ne0$ the case
$k=0$ is the round sphere $S^3(\tau^2)$.

Throughout, a \emph{surface} is a smooth, connected, oriented immersion of a
$2$-manifold $\Sigma$ into $\E$, and $\Sigma$ is \emph{closed} if it is compact
without boundary. Let $N$ be the unit normal,
$S=-\overline\nabla N$ the shape operator, $H=\tfrac12\operatorname{tr}S$ the mean
curvature, $u=\langle N,\xi\rangle$ the \emph{angle function}, and
$T=\xi-uN$ the tangent part of $\xi$. In a conformal parameter $\zeta$, put
$g=\lambda|d\zeta|^2$, $p=\langle S\partial_\zeta,\partial_\zeta\rangle$,
and $A=\langle\xi,\partial_\zeta\rangle$. The \emph{original}
Abresch--Rosenberg differential is
\begin{equation}\label{eq:QAR}
 \mathcal Q_{\mathrm{AR}}=\QAR\,d\zeta^2,\qquad \QAR=2(H+i\tau)p-kA^2;
\end{equation}
in the product spaces $\QAR=2Hp-\kappa A^2$.
Abresch and Rosenberg proved that $\mathcal Q_{\mathrm{AR}}$ is holomorphic
on every CMC surface~\cite[Theorem~1, p.~142]{AR2004}, \cite[Theorems~4 and~5, pp.~16--17]{AR2005}. Our main result is the converse, without any global
hypothesis.

\begin{theorem}\label{thm:main}
Let $\tau\ne0$ and $\kappa\in\R$. A surface in $\E$ has holomorphic
$\mathcal Q_{\mathrm{AR}}$ if and only if its mean curvature is constant.
\end{theorem}

No completeness, compactness, or analyticity is assumed, and the round case
$\kappa=4\tau^2$ is included; in particular, for $\tau\ne0$ there are no local
non-CMC examples. The hypothesis $\tau\ne0$ is essential: in every product
space $M^2(\kappa)\times\R$ there are non-CMC surfaces with holomorphic
$\calQAR$ (Proposition~\ref{prop:products}). For $\kappa\ne0$ we determine all
of them.

\begin{theorem}\label{thm:productsintro}
Let $\kappa\ne0$, and let $\Sigma$ be a surface in $M^2(\kappa)\times\R$ with
holomorphic $\calQAR$ and nonconstant mean curvature. Then $dH$ vanishes
nowhere on $\Sigma$, and:
\begin{enumerate}
\item[(a)] the angle function is a constant $u_0\notin\{0,\pm1\}$;
\item[(b)] $\nabla H$ is a principal direction with principal curvature $0$,
and $T$ is parallel to $\nabla H$;
\item[(c)] $|\nabla H|=\dfrac{|u_0|\,\bigl|4H^2+\kappa(1-u_0^2)\bigr|}{2\sqrt{1-u_0^2}}$;
\item[(d)] $\Sigma$ is locally congruent to one of the explicit constant angle
surfaces of Proposition~\ref{prop:constantangle}, with the same $u_0$; up to
local congruence it is determined by $\kappa$, $u_0$ and the sign of
$4H^2+\kappa(1-u_0^2)$.
\end{enumerate}
Items (a)--(c) and the nonvanishing of $dH$ hold on all of $\Sigma$; item (d) is a
local statement.
\end{theorem}

The surfaces in (d) are written down explicitly in
Proposition~\ref{prop:constantangle}, and Theorem~\ref{thm:products} identifies
$\Sigma$ with one of them. For $\kappa=\pm1$, those for which
$4H^2+\kappa(1-u_0^2)>0$ are, after a change of conformal parameter, the surfaces of Proposition~\ref{prop:products}; those with $r>0$, hence $H<0$, are the examples of \cite[Example~3.1]{FM2007}, and for $\kappa>0$ the family continues smoothly across the curve $H=0$; the new
content is that constant angle is necessary and that this family exhausts all
possibilities, with no completeness assumption. Since the mean curvature of a
closed surface has a critical point, Theorem~\ref{thm:productsintro} gives:

\begin{corollary}\label{cor:closedproducts}
Let $\kappa\ne0$. Every closed surface in $M^2(\kappa)\times\R$ with holomorphic
$\calQAR$ has constant mean curvature.
\end{corollary}

The implication CMC $\Rightarrow$ holomorphic is due to Abresch and Rosenberg
\cite[Theorem~1, p.~142]{AR2004}, \cite[Theorems~4 and~5, pp.~16--17]{AR2005}.
For the converse, in $\mathrm{Nil}_3$ Berdinsky and Taimanov
\cite[Proposition~3, p.~1014]{BT2005} proved that the normalized differential
$\widetilde p=\QAR/(2H+i)$ is holomorphic only on CMC surfaces, and Dorfmeister,
Inoguchi and Kobayashi \cite[Theorem~A.1, p.~441]{DIK2016} showed that a surface
with holomorphic original differential has constant mean curvature or is a Hopf
cylinder; we show in Section~\ref{sec:context} that the second alternative
reduces to the first. For the normalized differential
$Q_{\mathrm{FM}}=\QAR/(H+i\tau)$, Fern\'andez and Mira
\cite[Theorem~1.1, p.~282]{FM2007} obtained a characterization in
$\widetilde{PSL}_2(\R)$, $\mathrm H^2\times\R$ and the Berger spheres with
$0<\kappa\le8\tau^2$, with exceptional non-CMC surfaces when the base curvature
is negative; the cases $S^2\times\R$ and $\kappa>8\tau^2$ are left open in
\cite[\S1, p.~282]{FM2007}. To our knowledge, Theorem~\ref{thm:main} is the
first characterization of CMC surfaces by the original differential valid for
every $\kappa$ and every $\tau\ne0$, and it admits no exceptional surfaces.
Under compactness, Fern\'andez and Mira proved that a compact surface in $\E$
with holomorphic $Q_{\mathrm{FM}}$, and with nonzero Euler characteristic if
$\tau\ne0$, has constant mean curvature
\cite[Remark~(3), p.~282, and \S4, p.~289]{FM2007}, and they note, referring to \cite{CL2007}, that for $\tau=0$
the same holds for the original differential \cite[p.~289]{FM2007}; in the product spaces,
Cavalcante and Lira \cite[pp.~172--173]{CL2007} showed that on the set where
$dH\ne0$ a surface with holomorphic original differential has zero extrinsic
curvature, with $\nabla H$ a principal direction of principal curvature $0$, and observed
that non-CMC examples must be noncompact. These arguments use the topology of
the surface, whereas Theorem~\ref{thm:main} is local; by
Proposition~\ref{prop:products}, no local statement of this kind can hold when
$\tau=0$. Table~\ref{Table2} in Section~\ref{sec:context} compares these
results.

For closed surfaces, Theorem~\ref{thm:main} and the degree of the bundle of
quadratic differentials give:

\begin{corollary}\label{cor:zeros}
Let $\tau\ne0$, and let $\Sigma$ be a closed surface of genus $\mathsf g$ in $\E$
with holomorphic $\mathcal Q_{\mathrm{AR}}$. Then $H$ is constant, and:
\begin{enumerate}
\item[(i)] if $\mathsf g=0$, then $\mathcal Q_{\mathrm{AR}}\equiv0$ and $\Sigma$ is
a rotationally invariant CMC sphere;
\item[(ii)] if $\mathsf g=1$, then $\mathcal Q_{\mathrm{AR}}$ either vanishes
identically or has no zeros;
\item[(iii)] if $\mathsf g\ge2$ and $\mathcal Q_{\mathrm{AR}}\not\equiv0$, then
$\mathcal Q_{\mathrm{AR}}$ has exactly $4\mathsf g-4$ zeros, counted with
multiplicity.
\end{enumerate}
\end{corollary}

The rotational invariance in (i) follows from \cite[Theorem~6, p.~18]{AR2005}
(proof deferred there), \cite[Lemma~6.1 and Theorem~2.3]{ER2011} and
\cite[Proposition~2.4]{DVM2021}. We do not claim that these spheres are
embedded: in some Berger spheres they are only immersed
\cite[p.~594, and Corollary~1, p.~601]{Torralbo2010}.

\subsection{Strategy and organization.}
Section~\ref{sec:context} compares the normalizations of the
Abresch--Rosenberg differential used in the literature and the results closest
to ours. Section~\ref{sec:formalism} fixes the conventions, compares the
holomorphy conditions of $\QAR$ and $Q_{\mathrm{FM}}$
(Proposition~\ref{prop:comparison}), and performs the local reduction:
if $\mathcal Q_{\mathrm{AR}}$ is holomorphic and
$\Omega=\{dH\ne0\}\ne\emptyset$, then in the frame
$e_1=\nabla H/|\nabla H|$, $e_2=Je_1$,
\[
 [S]_{(e_1,e_2)}=\begin{pmatrix}0&\tau\\\tau&2H\end{pmatrix},\qquad \det S=-\tau^2,\qquad K=ku^2,
\]
where $K$ denotes the intrinsic curvature (Proposition~\ref{prop:normalform}).
Geometrically, $N$ is parallel along $\nabla H$ for a metric connection $D$
which preserves $\xi$ and whose curvature is $k$ times a rotation about $\xi$
(\S\ref{sec:connectionD}); $D$ is flat exactly in the round case $k=0$.
Section~\ref{sec:degenerate} settles the round case by an elementary
argument. Section~\ref{sec:puiseuxsetup} treats $k\ne0$: the structure
equations become a rational system whose solutions are algebraic, and an
embedding into a field of Puiseux series excludes every branch.
Section~\ref{sec:products} treats the product spaces, Section~\ref{sec:global}
proves Corollary~\ref{cor:zeros}, and the Appendix contains the computation of
the coefficients of the compatibility identity and facts on derivations.

\section{Normalizations and related results}\label{sec:context}

The Fern\'andez--Mira differential \cite[\S2, pp.~283--285, and Remark~2.6, p.~285]{FM2007} is
\begin{equation}\label{eq:QFM}
 Q_{\mathrm{FM}}=2p-\frac{k}{H+i\tau}A^2,\qquad \QAR=(H+i\tau)Q_{\mathrm{FM}}.
\end{equation}
Here $Q_{\mathrm{FM}}$ is defined where $H+i\tau\ne0$, that is, everywhere if
$\tau\ne0$ and on $\{H\ne0\}$ if $\tau=0$. The factor $H+i\tau$ is constant on the
CMC class, and nonzero there unless $\tau=H=0$; but off the CMC class
the two holomorphy conditions differ (Proposition~\ref{prop:comparison}).
In the round case $\QAR=2(H+i\tau)p$, not merely the Hopf coefficient $p$.
Table~\ref{Table1} relates the coefficients used in the literature.

\begin{table}[t]
\centering
\small
\renewcommand{\arraystretch}{1.3}
\begin{tabular}{|>{\raggedright\arraybackslash}p{3.2cm}|>{\centering\arraybackslash}p{4.4cm}|>{\centering\arraybackslash}p{3.0cm}|}
\hline
Reference & Coefficient used & Relation to $Q_{\mathrm{AR}}$\\
\hline
Abresch--Rosenberg \cite[p.~142]{AR2004}, \cite[pp.~16--17]{AR2005} & $Q_{\mathrm{AR}}=2(H{+}i\tau)p-kA^2$ & $Q_{\mathrm{AR}}$ itself\\
\hline
Fern\'andez--Mira \cite[pp.~284--285]{FM2007} & $Q_{\mathrm{FM}}=2p-\tfrac{k}{H+i\tau}A^2$ & $Q_{\mathrm{AR}}=(H{+}i\tau)Q_{\mathrm{FM}}$\\
\hline
Torralbo \cite[p.~594]{Torralbo2010} ($\kappa\ne4\tau^2$) & $\Theta=p-\tfrac{k}{2(H+i\tau)}A^2=\tfrac12Q_{\mathrm{FM}}$ & $Q_{\mathrm{AR}}=2(H{+}i\tau)\Theta$\\
\hline
Berdinsky--Taimanov \cite[p.~1014]{BT2005} ($\mathrm{Nil}_3$) & $\widetilde p=p+\tfrac{A^2}{2H+i}$ & $Q_{\mathrm{AR}}=(2H{+}i)\widetilde p$\\
\hline
Dorfmeister--Inoguchi--Kobayashi \cite[eq.~(3.9), p.~420]{DIK2016} ($\mathrm{Nil}_3$) & $B=\tfrac14(2H{+}i)\widetilde A$ & $Q_{\mathrm{AR}}=4B$\\
\hline
\end{tabular}
\caption{Coefficients used for the Abresch--Rosenberg differential in the literature, translated into the notation and sign conventions of the present paper, and their relation to the original differential $Q_{\mathrm{AR}}=2(H+i\tau)p-kA^2$ of \eqref{eq:QAR}. For $\mathrm{Nil}_3$ the normalization is $\kappa=0$, $\tau=\tfrac12$, so $k=-1$; $\widetilde A$ is the notation of \cite{DIK2016}.}\label{Table1}
\end{table}

For $\mathrm{Nil}_3$, Dorfmeister, Inoguchi and Kobayashi
\cite[Theorem~A.1, p.~441]{DIK2016} prove that if $B=\QAR/4$ is holomorphic,
the surface is CMC or a Hopf cylinder, and \cite[p.~411]{DIK2016} presents the
Hopf cylinder as a non-CMC example. For the original differential the second
alternative contains no non-CMC surface: on a Hopf cylinder $\QAR=H^2+\kappa/4$
depends only on the arclength of the base curve
(Remark~\ref{rem:hopfcylinder}), so holomorphy forces $H$ constant.
Table~\ref{Table2} summarizes the results closest to ours; by
Proposition~\ref{prop:comparison}, the holomorphy conditions of $\QAR$ and
$Q_{\mathrm{FM}}$ differ off the CMC class, and so do their exceptional
surfaces (Section~\ref{sec:products}).

\begin{table}[t]
\centering
\small
\renewcommand{\arraystretch}{1.3}
\begin{tabular}{|>{\raggedright\arraybackslash}p{1.9cm}|>{\raggedright\arraybackslash}p{2.5cm}|>{\centering\arraybackslash}p{1.75cm}|>{\centering\arraybackslash}p{1.25cm}|>{\raggedright\arraybackslash}p{2.55cm}|}
\hline
Result & Space & Differential & Setting & Conclusion\\
\hline
\cite{AR2004}, \cite{AR2005} & $\E$ & $\QAR$ & local & $H$ constant $\Rightarrow$ holomorphic\\
\hline
\cite[Prop.~3]{BT2005} & $\mathrm{Nil}_3$ & $\QAR/(2H+i)$ & local & holomorphic $\Rightarrow$ CMC\\
\hline
\cite[Thm.~A.1]{DIK2016} & $\mathrm{Nil}_3$ & $\QAR$ & local & holomorphic $\Rightarrow$ CMC or Hopf cylinder\\
\hline
\cite[Thm.~1.1]{FM2007} & $\widetilde{PSL}_2(\R)$, $\mathrm H^2\times\R$, Berger spheres with $0<\kappa\le8\tau^2$ & $Q_{\mathrm{FM}}$ & local & holomorphic $\Rightarrow$ CMC, or exceptional non-CMC surfaces if $\kappa<0$\\
\hline
\cite[Rem.~(3)]{FM2007} & $\E$ (with $\chi(\Sigma)\ne0$ if $\tau\ne0$) & $Q_{\mathrm{FM}}$ & compact & holomorphic $\Rightarrow$ CMC\\
\hline
\cite[p.~289]{FM2007}, \cite{CL2007} & $M^2(\kappa)\times\R$ & $\QAR$ & compact & holomorphic $\Rightarrow$ CMC\\
\hline
Theorem~\ref{thm:main} & $\E$, $\tau\ne0$, all $\kappa$ & $\QAR$ & local & holomorphic $\Leftrightarrow$ CMC\\
\hline
Theorem~\ref{thm:productsintro} & $M^2(\kappa)\times\R$, $\kappa\ne0$ & $\QAR$ & local & non-CMC $\Rightarrow$ constant angle, explicit family\\
\hline
\end{tabular}
\caption{Results closest to the present paper. ``Local'' means that no completeness or compactness is assumed.}\label{Table2}
\end{table}

Daniel's sister correspondence \cite[Theorem~5.2, p.~114]{Daniel2007}
multiplies $\QAR$ by a unimodular constant but is defined only for CMC surfaces,
so it gives no information on the converse. For a survey see \cite{FM2010}; for
a Codazzi pair whose $(2,0)$-part is $\calQAR$, see \cite{ET2019}.

The companion paper \cite[Theorem~1.1, p.~2]{AlencarRosenberg2026b} shows that the
inequality $|dH|_g\le h|\mathcal Q_{\mathrm{AR}}|_g$, with $h\ge0$ continuous (or
in $L^p_{\mathrm{loc}}$ for some $p>2$), forces $\calQAR\equiv0$ and $H$ constant
on topological spheres in every $\E$, the round case $\kappa=4\tau^2$ included;
for $\tau=0$ and for $\tau\ne0$, $\kappa\ne4\tau^2$, this conclusion is due to
\cite[Theorem~1.1, p.~284]{ACT2007} and \cite[Theorem~1, p.~526]{ACFT2007}. That
sphere theorem uses the topology of the sphere and does not yield the local
converse proved here. Moreover, the genus-zero hypothesis cannot be replaced by
any other genus: in the round sphere and in every Berger sphere there are compact
embedded non-CMC surfaces of every genus $g\ge1$ satisfying the inequality with a
constant $h$ \cite[Corollary~6.4, p.~15]{AlencarRosenberg2026b}.

\section{Geometric formalism and local reduction}\label{sec:formalism}

\emph{Conventions.} We retain the notation of the introduction, including
$k=\kappa-4\tau^2$ and the round interpretation when $k=0$. The orientation of
$\E$ is the one for which $\overline\nabla_X\xi=\tau X\times\xi$; the
rotation of the tangent plane is $JX=N\times X$; conformal parameters
$\zeta=x+iy$ are positively oriented, with
$\partial_\zeta=\tfrac12(\partial_x-i\partial_y)$; for any connection the
curvature is $R(X,Y)=\nabla_X\nabla_Y-\nabla_Y\nabla_X-\nabla_{[X,Y]}$; and $K$
is the intrinsic (Gaussian) curvature of the induced metric. Two surfaces are
\emph{locally congruent} if every point of the first has a neighborhood that an
isometry of the ambient space carries onto an open subset of the second. The
behavior of these conventions under homotheties and under the change of
orientation $\tau\mapsto-\tau$ is recorded before \eqref{eq:normalization}.

Daniel's structure equations \cite[Corollary~3.2 and Proposition~3.3, pp.~97--98; Theorem~4.3, p.~101]{Daniel2007} are
\begin{align}
 K&=\det S+\tau^2+ku^2,\label{eq:gauss}\\
 (\nabla_XS)Y-(\nabla_YS)X&=ku(\langle Y,T\rangle X-\langle X,T\rangle Y),\label{eq:codazzi}\\
 \nabla_XT&=u(SX-\tau JX),\label{eq:T}\\
 X(u)&=-\langle SX-\tau JX,T\rangle,\qquad |T|^2+u^2=1.\label{eq:vertical}
\end{align}
Here $K$ is the intrinsic curvature. In conformal coordinates,
\begin{equation}\label{eq:conformal}
 K=-\frac{\Delta_0\log\lambda}{2\lambda},\qquad \det S=H^2-\frac{4|p|^2}{\lambda^2},\qquad \Delta_0=4\partial_\zeta\partial_{\bar\zeta}.
\end{equation}
The first identity is Liouville's curvature formula for isothermal
parameters, see \cite[\S4-3, Exercise~1]{doCarmo1976}. For the second, with
$e_1=\partial_x/\sqrt\lambda$, $e_2=\partial_y/\sqrt\lambda$ and
$S_{ij}=\langle Se_i,e_j\rangle$, bilinearity and
$\partial_\zeta=\tfrac{\sqrt\lambda}2(e_1-ie_2)$ give
$p=\tfrac\lambda4(S_{11}-S_{22}-2iS_{12})$, hence
$4|p|^2/\lambda^2=\tfrac14(S_{11}-S_{22})^2+S_{12}^2=H^2-\det S$.
In particular, $p=0$ precisely at umbilic points.

The fundamental equations \cite[Theorem~2.3, p.~283]{FM2007} are
\begin{align}
 A_{\bar\zeta}&=\frac{u\lambda}{2}(H+i\tau), & p_{\bar\zeta}&=\frac\lambda2(H_\zeta+kuA),\label{eq:fund1}\\
 u_\zeta&=-(H-i\tau)A-\frac{2p}{\lambda}\overline A, & 4|A|^2&=\lambda(1-u^2).\label{eq:fund2}
\end{align}
They follow from \eqref{eq:codazzi}--\eqref{eq:vertical} by complex bilinear
extension, using $J\partial_\zeta=i\partial_\zeta$,
$\langle S\partial_\zeta,\partial_{\bar\zeta}\rangle=\lambda H/2$,
$T=\tfrac{2\bar A}\lambda\partial_\zeta+\tfrac{2A}\lambda\partial_{\bar\zeta}$ and
the Christoffel symbols of the conformal metric,
\begin{equation}\label{eq:christoffelconformal}
 \nabla_{\partial_\zeta}\partial_\zeta=(\log\lambda)_\zeta\,\partial_\zeta,\qquad \nabla_{\partial_\zeta}\partial_{\bar\zeta}=\nabla_{\partial_{\bar\zeta}}\partial_\zeta=0.
\end{equation}
For instance, $A_{\bar\zeta}=\langle\nabla_{\partial_{\bar\zeta}}T,\partial_\zeta\rangle$
with $\nabla_{\partial_{\bar\zeta}}T=u(S\partial_{\bar\zeta}+i\tau\partial_{\bar\zeta})$
by \eqref{eq:T}, and $p_{\bar\zeta}$ is obtained by pairing the complexified
Codazzi equation $(\nabla_\zeta S)\partial_{\bar\zeta}-(\nabla_{\bar\zeta}S)\partial_\zeta
=ku(\bar A\partial_\zeta-A\partial_{\bar\zeta})$ with $\partial_\zeta$. The signs
agree with \cite[(C.1)--(C.4)]{FM2007}.

Differentiating \eqref{eq:QAR} and using \eqref{eq:fund1} gives
\begin{align*}
 (\QAR)_{\bar\zeta}&=2H_{\bar\zeta}\,p+2(H+i\tau)p_{\bar\zeta}-2kA\,A_{\bar\zeta}\\
 &=2H_{\bar\zeta}p+\lambda(H+i\tau)H_\zeta+ku\lambda A(H+i\tau)-ku\lambda A(H+i\tau).
\end{align*}
The terms $ku\lambda A(H+i\tau)$ cancel, leaving, for every immersion,
\begin{equation}\label{eq:holoidentity}
 (\QAR)_{\bar\zeta}=2pH_{\bar\zeta}+\lambda(H+i\tau)H_\zeta.
\end{equation}
Thus holomorphy is equivalent to
\begin{equation}\label{eq:holo}
 2pH_{\bar\zeta}+\lambda(H+i\tau)H_\zeta=0.
\end{equation}
If $H$ is constant, both derivatives vanish, proving the CMC-to-holomorphy
implication in Theorem~\ref{thm:main}. The identity \eqref{eq:holoidentity},
valid for every immersion, is the starting point for the converse: on the
open set where $dH\ne0$, equation \eqref{eq:holo} determines the Hopf
coefficient $p$ in terms of $H$ and $\lambda$, and this rigidly fixes the
shape operator (Proposition~\ref{prop:normalform} below).

\begin{proposition}[Comparison of the holomorphy conditions]\label{prop:comparison}
For any immersion with $\tau\ne0$, with no CMC hypothesis,
\[
 (Q_{\mathrm{FM}})_{\bar\zeta}=\lambda H_\zeta+\frac{kA^2}{(H+i\tau)^2}H_{\bar\zeta},\qquad
 (\QAR)_{\bar\zeta}=(H+i\tau)(Q_{\mathrm{FM}})_{\bar\zeta}+H_{\bar\zeta}\,Q_{\mathrm{FM}}.
\]
Consequently the holomorphy equations of $\QAR$ and $Q_{\mathrm{FM}}$ are, respectively,
\[
 2pH_{\bar\zeta}+\lambda(H+i\tau)H_\zeta=0,\qquad kA^2H_{\bar\zeta}+\lambda(H+i\tau)^2H_\zeta=0,
\]
and $\QAR$, $Q_{\mathrm{FM}}$ have the same zero set.
\end{proposition}
\begin{proof}
Differentiating \eqref{eq:QFM},
\[
 (Q_{\mathrm{FM}})_{\bar\zeta}=2p_{\bar\zeta}-k\left[\frac{2AA_{\bar\zeta}}{H+i\tau}-\frac{A^2H_{\bar\zeta}}{(H+i\tau)^2}\right].
\]
By \eqref{eq:fund1}, $2p_{\bar\zeta}=\lambda H_\zeta+ku\lambda A$ and
\[
 \frac{2kAA_{\bar\zeta}}{H+i\tau}=\frac{2kA}{H+i\tau}\cdot\frac{u\lambda}2(H+i\tau)=ku\lambda A .
\]
The two terms $ku\lambda A$ cancel, which gives the first identity. The second is the Leibniz rule applied to $\QAR=(H+i\tau)Q_{\mathrm{FM}}$, using $(H+i\tau)_{\bar\zeta}=H_{\bar\zeta}$. The holomorphy equation of $\QAR$ is \eqref{eq:holo}; that of $Q_{\mathrm{FM}}$ follows from the first identity after multiplication by $(H+i\tau)^2\ne0$ ($H$ real, $\tau\ne0$). Finally $H+i\tau\ne0$ everywhere, so $\QAR=0\iff Q_{\mathrm{FM}}=0$.
\end{proof}

Equality of the zero sets does not make the two holomorphy conditions
equivalent: if $Q_{\mathrm{FM}}$ is holomorphic, then
$(\QAR)_{\bar\zeta}=H_{\bar\zeta}Q_{\mathrm{FM}}$. Fern\'andez--Mira
\cite[Remark~2.6, p.~285]{FM2007} point out this distinction (they write
$P\,d\zeta^2=(H+i\tau)Q\,d\zeta^2$ for the original differential); by
Theorem~\ref{thm:main}, their exceptional surfaces in $\widetilde{PSL}_2(\R)$
\cite[Theorem~1.1]{FM2007} do not have holomorphic $\QAR$.

\begin{remark}[Hopf cylinders with holomorphic $\calQAR$ are CMC]\label{rem:hopfcylinder}
Independently of Theorem~\ref{thm:main}: a Hopf cylinder $\pi^{-1}(\gamma)$ is
characterized by $u\equiv0$, so $T=\xi$ is parallel by \eqref{eq:T}; with the
parallel frame $e_2=\xi$, $e_1=-Je_2$ and flat coordinates
($\partial_x=e_1$, $\partial_y=e_2$, $\lambda=1$), \eqref{eq:vertical} gives
$S_{22}=0$, $S_{12}=\tau$, hence $p=(H-i\tau)/2$, $A=-i/2$ and
\[
 \QAR=(H+i\tau)(H-i\tau)+\frac k4=H^2+\frac\kappa4 .
\]
Since $H$ is constant along the fibers, every datum is a function of $x$ alone,
and a holomorphic function of $x$ alone is constant; so holomorphy forces $H^2$,
hence $H$, to be constant. The same argument applies to the normalized
differentials, e.g.\ $\QAR/(2H+i)=H^2/(2H+i)$ in $\mathrm{Nil}_3$.
\end{remark}

Assume from now on that $\mathcal Q_{\mathrm{AR}}$ is holomorphic and that $\Omega=\{dH\ne0\}$ is nonempty; write $f=|\nabla H|>0$, $e_1=\nabla H/f$, $e_2=Je_1$, a globally defined smooth positively oriented orthonormal frame on a fixed connected component of $\Omega$, which we continue to denote by $\Omega$.

\begin{proposition}[Normal form on the non-CMC set]\label{prop:normalform}
If $\mathcal Q_{\mathrm{AR}}$ is holomorphic, then on $\Omega$
\begin{equation}\label{eq:normalform}
 [S]_{(e_1,e_2)}=\begin{pmatrix}0&\tau\\\tau&2H\end{pmatrix},\qquad \det S=-\tau^2,\qquad K=ku^2.
\end{equation}
In particular, $\nabla H$ is an asymptotic direction at every point of
$\Omega$, and the principal curvatures are $H\pm\sqrt{H^2+\tau^2}$.
\end{proposition}
\begin{proof}
At $x_0\in\Omega$, rotate the conformal parameter, $\zeta\mapsto e^{i\theta_0}\zeta$,
so that $H_\zeta(x_0)>0$; then $\lambda$ is unchanged, \eqref{eq:holo} is
invariant (both terms acquire the factor $e^{-i\theta_0}$), and at $x_0$ the
frame $(\partial_x,\partial_y)/\sqrt\lambda$ is $(e_1,e_2)$, so
$H_\zeta=H_{\bar\zeta}=f\sqrt\lambda/2$. Then \eqref{eq:holo} gives
$p(x_0)=-\lambda(x_0)(H(x_0)+i\tau)/2$. Since
$p=\tfrac\lambda4(S_{11}-S_{22}-2iS_{12})$ (from
$\partial_\zeta=\tfrac{\sqrt\lambda}2(e_1-ie_2)$), this gives $S_{11}-S_{22}=-2H$
and $S_{12}=\tau$; with $S_{11}+S_{22}=2H$, $S_{11}=0$ and $S_{22}=2H$.
Then $\det S=-\tau^2$, and \eqref{eq:gauss} gives $K=\det S+\tau^2+ku^2=ku^2$.
Finally $\langle Se_1,e_1\rangle=S_{11}=0$, so $e_1=\nabla H/|\nabla H|$ is
asymptotic, and the eigenvalues of $S$ are
$H\pm\sqrt{H^2-\det S}=H\pm\sqrt{H^2+\tau^2}$.
\end{proof}

The normal form is a pointwise consequence of \eqref{eq:holo} and
\eqref{eq:gauss}; no analyticity or global hypothesis is used. Thus every point
of $\Omega$ is hyperbolic, with $\kappa_1\kappa_2=-\tau^2$, and $\nabla H$ is
asymptotic but not principal ($\langle Se_1,e_2\rangle=\tau\ne0$), in contrast
with the product case $\tau=0$, where $\nabla H$ is principal with principal
curvature $0$ \cite[pp.~172--173]{CL2007}. We now show that the structure
equations are incompatible with this normal form.

On $\Omega$ write $\nabla_{e_1}e_1=a\,e_2$ and $\nabla_{e_2}e_1=b\,e_2$, so that
$a,b$ are the coefficients of the Levi-Civita connection in the frame
$(e_1,e_2)$, and put $\rho^2=H^2+\tau^2>0$. On the open subset of $\Omega$ where
$R:=\sqrt{1-u^2}>0$, choose locally a smooth angle function $\theta$ with
\[
 T=\xi-uN=R(\cos\theta\,e_1+\sin\theta\,e_2),
\]
and put
\[
 P=\tau\cos\theta+H\sin\theta,\qquad L=-\tau\sin\theta+H\cos\theta,
\]
so that $P^2+L^2=\rho^2$.

\begin{lemma}[Frame equations]\label{lem:normal}
If $\mathcal Q_{\mathrm{AR}}$ is holomorphic, then $e_1u=0$ on $\Omega$, and,
where $R>0$,
\begin{align}
 a&=-\frac{\tau f}{\rho^2}-\frac{kuRP}{2\rho^2}, & b&=-\frac{Hf}{\rho^2}-\frac{kuRL}{2\rho^2},\label{eq:ab}\\
 e_1u&=0, & e_2u&=-2RP,\label{eq:eu}\\
 e_1\theta&=-a, & e_2\theta&=-b+\frac{2uL}{R},\qquad e_2f=af.\label{eq:etheta}
\end{align}
\end{lemma}
\begin{proof}
We use the normal form \eqref{eq:normalform} throughout.

\emph{Step 1: $e_1u=0$ and $e_2u=-2RP$.} Since $Se_1=\tau e_2$, $Se_2=\tau e_1+2He_2$,
$Je_1=e_2$ and $Je_2=-e_1$,
\begin{equation}\label{eq:StauJ}
 (S-\tau J)e_1=0,\qquad (S-\tau J)e_2=2(\tau e_1+He_2).
\end{equation}
By \eqref{eq:vertical}, $e_1u=-\langle(S-\tau J)e_1,T\rangle=0$. Where $R>0$, write
$T=R(\cos\theta\,e_1+\sin\theta\,e_2)$; then
$e_2u=-\langle2\tau e_1+2He_2,T\rangle=-2R(\tau\cos\theta+H\sin\theta)=-2RP$.

\emph{Step 2: the connection coefficients.} Write $\nabla_{e_1}e_1=ae_2$, $\nabla_{e_1}e_2=-ae_1$,
$\nabla_{e_2}e_1=be_2$, $\nabla_{e_2}e_2=-be_1$.
Using $(\nabla_XS)Y=\nabla_X(SY)-S\nabla_XY$ and
$e_1H=f$, $e_2H=0$, we obtain
\begin{align*}
 (\nabla_{e_1}S)e_2&=\nabla_{e_1}(\tau e_1+2He_2)-S(-ae_1)\\
 &=\tau ae_2+2fe_2-2Hae_1+a\tau e_2=-2Ha\,e_1+(2\tau a+2f)e_2,\\
 (\nabla_{e_2}S)e_1&=\nabla_{e_2}(\tau e_2)-S(be_2)\\
 &=-\tau b\,e_1-b(\tau e_1+2He_2)=-2\tau b\,e_1-2Hb\,e_2.
\end{align*}
Subtracting and comparing with the right side
$kuR\sin\theta\,e_1-kuR\cos\theta\,e_2$ of \eqref{eq:codazzi} gives
$-2Ha+2\tau b=kuR\sin\theta$ and $2\tau a+2Hb=-2f-kuR\cos\theta$, a linear
system of determinant $-4\rho^2\ne0$ whose solution is \eqref{eq:ab}.

\emph{Step 3: the angular equations.} Since $e_1u=0$, $e_1R=0$, and
$\nabla_{e_1}T=R(e_1\theta+a)(-\sin\theta\,e_1+\cos\theta\,e_2)$; by \eqref{eq:T}
and \eqref{eq:StauJ}, $\nabla_{e_1}T=u(S-\tau J)e_1=0$, so $e_1\theta=-a$.
Similarly $e_2R=2uP$, and the component of $\nabla_{e_2}T$ along
$-\sin\theta\,e_1+\cos\theta\,e_2$ is $R(e_2\theta+b)$, while
$u(S-\tau J)e_2=2u(\tau e_1+He_2)$ has component $2uL$; hence
$e_2\theta=-b+2uL/R$.

\emph{Step 4: $e_2f=af$.} Since $f=e_1H$, $e_2f=e_2e_1H$. By torsion-freeness, $[e_2,e_1]=\nabla_{e_2}e_1-\nabla_{e_1}e_2=be_2+ae_1$, so $e_2e_1H-e_1e_2H=[e_2,e_1]H=af$. Since $e_2H\equiv0$ on $\Omega$, $e_1e_2H=0$, and hence $e_2f=e_2e_1H=af$.
\end{proof}

Put $\Omega':=\Omega\cap\{|u|<1\}$, where $R>0$, and
\[
 \Omega_*:=\{x\in\Omega':\ H(x)\ne0,\ u(x)\ne0,\ P(x)\ne0\}.
\]
On $\Omega'$ the unit vector $T/R$, and with it $\cos\theta$, $\sin\theta$, $P$
and $L$, are well defined and smooth, so $\Omega_*$ is open.

\begin{lemma}[Regular coordinate region]\label{lem:density}
If $\mathcal Q_{\mathrm{AR}}$ is holomorphic, then $\Omega_*$ is dense in
$\Omega$, and $(H,u)$ are local coordinates on $\Omega_*$.
\end{lemma}
\begin{proof}
Since the complement of a closed set with empty interior is open and dense,
it suffices to show that $\{|u|=1\}$ has empty interior in $\Omega$ and that
the relatively closed sets $\{H=0\}$, $\{u=0\}$ and $\{P=0\}$ of $\Omega'$ have
empty interior in $\Omega'$.

\emph{Step 1: $\{|u|=1\}$.} If $|u|\equiv1$ on a nonempty open
$U\subset\Omega$, then $|T|^2=1-u^2\equiv0$, so $\nabla T\equiv0$ on $U$. By
\eqref{eq:T} and \eqref{eq:StauJ}, $\nabla_{e_2}T=2u(\tau e_1+He_2)\ne0$,
since $u=\pm1$ and $\tau\ne0$, a contradiction.

\emph{Step 2: $\{H=0\}$.} If $H\equiv0$ on an open $U\subset\Omega'$, then
$dH=0$ there, contradicting $U\subset\Omega$.

\emph{Step 3: $\{P=0\}$.} Suppose $P\equiv0$ on a nonempty open
$U\subset\Omega'$; by Step~2 we may shrink $U$ so that $H\ne0$ on $U$. We derive
a contradiction in three substeps.

(i) \emph{$u$ and $\theta$ are constant along $e_2$.} By \eqref{eq:eu},
$e_2u=-2RP=0$. Since $e_2H=0$,
\[
 0=e_2P=(-\tau\sin\theta+H\cos\theta)\,e_2\theta=L\,e_2\theta ,
\]
and $L\ne0$ because $P^2+L^2=\rho^2>0$; hence $e_2\theta=0$.

(ii) \emph{$f$ is a function of $H$, $u$ and $\theta$.} By (i) and
\eqref{eq:etheta}, $b=2uL/R$. Comparing with \eqref{eq:ab} and solving for $f$,
\[
 \frac{2uL}R=-\frac{Hf}{\rho^2}-\frac{kuRL}{2\rho^2}\ \Longrightarrow\ f=-\frac{uL(4\rho^2+kR^2)}{2HR}.
\]

(iii) \emph{Contradiction.} Every factor on the right of (ii) has vanishing
$e_2$-derivative: $e_2u=0$, $e_2\rho^2=2H\,e_2H=0$, $e_2R=-u\,e_2u/R=0$, and
$e_2L=(e_2H)\cos\theta-P\,e_2\theta=0$. Hence $e_2f=0$. On the other hand,
\eqref{eq:etheta} and \eqref{eq:ab} with $P=0$ give
$e_2f=af=-\tau f^2/\rho^2\ne0$.

\emph{Step 4: $\{u=0\}$.} If $u\equiv0$ on a nonempty open $U\subset\Omega'$,
then $R=1$ and \eqref{eq:eu} gives $P\equiv0$ on $U$, which is excluded by
Step~3.

Hence $\Omega_*$ is dense in $\Omega$. Finally,
$(dH\wedge du)(e_1,e_2)=f\cdot(-2RP)-0=-2fRP\ne0$ on $\Omega_*$, and the inverse
function theorem shows that $(H,u)$ are local coordinates there.
\end{proof}

Under the homothety $\bar g\mapsto \mathfrak c^2\bar g$, $\mathfrak c>0$,
\[
(\kappa,\tau,H,\lambda,p,A)
\longmapsto(\kappa/\mathfrak c^2,\tau/\mathfrak c,H/\mathfrak c,\mathfrak c^2\lambda,\mathfrak c p,\mathfrak c A).
\]
Indeed, $S\mapsto S/\mathfrak c$, so $p\mapsto \mathfrak c^2p/\mathfrak c=\mathfrak c p$;
the vertical component scales by $\mathfrak c$, while $u$ is unchanged.
The ambient sectional curvatures scale by $\mathfrak c^{-2}$. Consequently,
\[
 2\Bigl(\frac H{\mathfrak c}+i\frac\tau{\mathfrak c}\Bigr)(\mathfrak c p)-\Bigl(\frac{\kappa-4\tau^2}{\mathfrak c^2}\Bigr)(\mathfrak c A)^2=2(H+i\tau)p-kA^2=\QAR.
\]
Take $\mathfrak c=|\tau|$. If $\tau<0$, recall that reversing the
orientation of $E(\kappa,\tau)$ gives an orientation-reversing isometry
$E(\kappa,\tau)\to E(\kappa,-\tau)$ (the sign of $\tau$ in
$\overline\nabla_X\xi=\tau X\times\xi$ depends on the orientation through
the cross product). Reverse accordingly the ambient and induced surface
orientations, keeping $N,\xi$. This reverses $\tau$,
replaces $\zeta$ by $\bar\zeta$, and conjugates $p,A,\QAR$,
preserving holomorphy. For $k\ne0$ we may therefore assume
\begin{equation}\label{eq:normalization}
 \tau=1,\qquad q=\frac{\kappa-4\tau^2}{\tau^2}\ne0.
\end{equation}
The parameter $q=k/\tau^2$ is invariant under both operations.

\subsection{A connection adapted to the fibration}\label{sec:connectionD}
We close this section with a geometric reading of the normal form, which
explains the role of $k$ and the dichotomy between Sections~\ref{sec:degenerate}
and~\ref{sec:puiseuxsetup}. It is not used in the proof of
Theorem~\ref{thm:main}. Orient $\E$ so that
$\overline\nabla_X\xi=\tau X\times\xi$, and define
\begin{equation}\label{eq:Ddef}
 D_XY:=\overline\nabla_XY-\tau\,X\times Y .
\end{equation}
Since $Y\mapsto X\times Y$ is skew-adjoint, $D$ is a metric connection; its
torsion $-2\tau X\times Y$ is totally skew-symmetric, and
$D_X\xi=\tau X\times\xi-\tau X\times\xi=0$. Thus $\xi$ is $D$-parallel, and
$D$ preserves the horizontal distribution $\xi^\perp$. Moreover $D$ commutes
with $Z\mapsto\xi\times Z$: expanding $D_X(\xi\times Y)-\xi\times D_XY$ with
$\overline\nabla_X\xi=\tau X\times\xi$ leaves
$\tau[(X\times\xi)\times Y+\xi\times(X\times Y)+X\times(Y\times\xi)]$, which
vanishes by the Jacobi identity.

\begin{proposition}[Curvature of $D$]\label{prop:curvD}
The curvature of $D$ is
\[
 R^D(X,Y)Z=-k\,\langle X\times Y,\xi\rangle\,\xi\times Z .
\]
In particular $D$ is flat if and only if $k=0$, and its holonomy consists of
rotations about $\xi$.
\end{proposition}
\begin{proof}
Write $D=\overline\nabla+\mathcal T$ with $\mathcal T_XY=-\tau X\times Y$. The
cross product is $\overline\nabla$-parallel, so $\overline\nabla\mathcal T=0$ and
$R^D(X,Y)=\overline R(X,Y)+[\mathcal T_X,\mathcal T_Y]$. By the Jacobi
identity, $[\mathcal T_X,\mathcal T_Y]Z=\tau^2(X\times Y)\times Z
=-\tau^2(\langle Y,Z\rangle X-\langle X,Z\rangle Y)$. Translated to the present sign
convention, the curvature of $\E$ is \cite[Proposition~2.1]{Daniel2007}
\[
\begin{aligned}
 \overline R(X,Y)Z={}&(\kappa-3\tau^2)(\langle Y,Z\rangle X-\langle X,Z\rangle Y)\\
 &-k\bigl(\langle Y,\xi\rangle\langle Z,\xi\rangle X-\langle X,\xi\rangle\langle Z,\xi\rangle Y\\
 &\qquad+\langle X,\xi\rangle\langle Y,Z\rangle\xi-\langle Y,\xi\rangle\langle X,Z\rangle\xi\bigr).
\end{aligned}
\]
Hence $R^D(X,Y)Z=k\,\mathcal C(X,Y)Z$, where $\mathcal C(X,Y)Z$ is
$\langle Y,Z\rangle X-\langle X,Z\rangle Y$ minus the bracket above. A direct
check shows $\mathcal C(X,Y)Z=0$ whenever one of $X,Y,Z$ equals $\xi$, so
$\mathcal C$ only sees horizontal components; on a positively oriented
horizontal orthonormal pair $E_1,E_2$ (so $E_1\times E_2=\xi$,
$\xi\times E_1=E_2$, $\xi\times E_2=-E_1$) one has
$\mathcal C(E_1,E_2)E_1=-E_2=-\xi\times E_1$ and
$\mathcal C(E_1,E_2)E_2=E_1=-\xi\times E_2$. Both sides of the claimed formula
are trilinear, which proves it.
\end{proof}

For any immersion, $JX=N\times X$ gives
\begin{equation}\label{eq:DN}
 D_XN=-(S-\tau J)X,\qquad \det(S-\tau J)=\det S+\tau^2=K-ku^2 .
\end{equation}

\begin{proposition}\label{prop:Dparallel}
On $\Omega=\{dH\ne0\}$, $\calQAR$ is holomorphic if and only if $D_{\nabla H}N=0$.
\end{proposition}
\begin{proof}
By \eqref{eq:DN}, $D_{e_1}N=0$ if and only if $Se_1=\tau e_2$, i.e.\ $S_{11}=0$,
$S_{12}=\tau$, which is the normal form \eqref{eq:normalform}; in the rotated
parameter of the proof of Proposition~\ref{prop:normalform} this is
$p=-\lambda(H+i\tau)/2$, i.e.\ \eqref{eq:holo}.
\end{proof}

On $\Omega'=\Omega\cap\{|u|<1\}$ let $\nu=(N-u\xi)/R$ be the normalized
horizontal part of $N$ and $\sigma(X)=\langle D_X\nu,\xi\times\nu\rangle$. Since $D$
preserves $\xi^\perp$ and commutes with $\xi\times$, Proposition~\ref{prop:curvD}
gives $d\sigma(X,Y)=-k\langle X\times Y,\xi\rangle$, and with
$\omega=a\,e_1^\flat+b\,e_2^\flat$ the connection form of $(e_1,e_2)$,
Lemma~\ref{lem:normal} reads
\begin{equation}\label{eq:dsigma}
 d\sigma=-ku\,dA,\qquad du=-2RP\,e_2^\flat,\qquad \sigma=\frac{2L}R\,e_2^\flat,\qquad d\theta=u\sigma-\omega .
\end{equation}
The integrability condition $d(d\theta)=0$ is the Gauss equation $K=ku^2$. For
$k=0$, $\sigma$ is closed, which is why the round case integrates directly
(Section~\ref{sec:degenerate}). For $k\ne0$ it is not; the curvature of $D$
alone does not exclude non-CMC surfaces, and the obstruction comes from the
full system of Section~\ref{sec:puiseuxsetup}.

\section{Proof of Theorem~\ref{thm:main} for \texorpdfstring{$k=0$}{k = 0}}\label{sec:degenerate}

In the round case $k=0$ the original differential is $\QAR=2(H+i\tau)p$,
and the algebraic reduction of Section~\ref{sec:puiseuxsetup} is not needed:
the proof uses only \eqref{eq:holo}, the Gauss equation and Liouville's
formula.

\begin{proof}[Proof of Theorem~\ref{thm:main} for $k=0$]
Suppose $\calQAR$ is holomorphic and, for contradiction, that $\Omega=\{dH\ne0\}\ne\emptyset$. Fix a point of $\Omega$ and a simply connected coordinate disc $\mathcal D\subset\Omega$ around it, with conformal parameter $\zeta=x+iy$, $g=\lambda|d\zeta|^2$. Put $\rho^2:=H^2+\tau^2>0$.

\emph{Step 1: a natural parameter.} Taking moduli in \eqref{eq:holo},
$|p|=\lambda\rho/2>0$, so $\QAR=2(H+i\tau)p$ has no zeros on $\mathcal D$. Choose a
holomorphic $\omega$ with $e^\omega=\QAR$ \cite[Ch.~4]{Ahlfors1966} and the new
conformal parameter $\eta=\int e^{\omega/2}d\zeta$, so that $d\eta^2=\calQAR$;
the relation $|p|=\lambda\rho/2$ is invariant, since $|p|$ and $\lambda$ acquire
the same factor. In this parameter $\QAR\equiv1$, so
\begin{equation}\label{eq:degplambda}
 p=\frac1{2(H+i\tau)},\qquad \lambda=\frac{2|p|}{\rho}=\frac1{\rho^2}.
\end{equation}

\emph{Step 2: the equation for $H$.} Substituting \eqref{eq:degplambda} in \eqref{eq:holo},
\[
 \frac{H_{\bar\zeta}}{H+i\tau}+\frac{(H+i\tau)H_\zeta}{\rho^2}=0 .
\]
Multiplying by $\rho^2=(H+i\tau)(H-i\tau)$ gives $(H-i\tau)H_{\bar\zeta}+(H+i\tau)H_\zeta=0$. As $H_{\bar\zeta}=\overline{H_\zeta}$, this says $2\operatorname{Re}\bigl[(H+i\tau)H_\zeta\bigr]=0$, and with $H_\zeta=\tfrac12(H_x-iH_y)$,
\begin{equation}\label{eq:degenerate}
 HH_x+\tau H_y=0 .
\end{equation}
The Codazzi equation $p_{\bar\zeta}=\lambda H_\zeta/2$ gives the same equation.

\emph{Step 3: Gauss.} By \eqref{eq:conformal} and \eqref{eq:degplambda}, $\det S=H^2-4|p|^2/\lambda^2=H^2-\rho^2=-\tau^2$, so \eqref{eq:gauss} gives $K=0$. By Liouville's formula $\Delta_0\log\lambda=0$, where $\Delta_0=\partial_x^2+\partial_y^2$ is the coordinate Laplacian of \eqref{eq:conformal}; since $\log\lambda=-\log\rho^2$,
\begin{equation}\label{eq:degenerate2}
 \Delta_0\psi=0,\qquad \psi:=\log(H^2+\tau^2).
\end{equation}

\emph{Step 4: computation.} Put $\beta:=H_x$. By \eqref{eq:degenerate}, $H_y=-H\beta/\tau$. Differentiating \eqref{eq:degenerate} in $x$,
\begin{equation}\label{eq:betay}
 \tau\beta_y+\beta^2+H\beta_x=0 .
\end{equation}
Differentiating $H_y=-H\beta/\tau$ in $y$ and using \eqref{eq:betay},
\[
 H_{yy}=-\frac{H_y\beta+H\beta_y}\tau=\frac{H\beta^2}{\tau^2}-\frac{H\beta_y}\tau=\frac H{\tau^2}\bigl(2\beta^2+H\beta_x\bigr).
\]
Moreover, for the coordinate gradient $\nabla_0H:=(H_x,H_y)$, $|\nabla_0H|^2=\beta^2+H^2\beta^2/\tau^2=\beta^2\rho^2/\tau^2$. Since $\psi$ is a function of $H$ with $\psi'=2H/\rho^2$ and $\psi''=2(\tau^2-H^2)/\rho^4$, the chain rule $\Delta_0\psi=\psi''|\nabla_0H|^2+\psi'\,\Delta_0H$, $\Delta_0H=\beta_x+H_{yy}$, gives
\[
 \begin{aligned}
 \Delta_0\psi&=\frac{2(\tau^2-H^2)\beta^2}{\tau^2\rho^2}+\frac{2H}{\rho^2}\Bigl[\beta_x+\frac H{\tau^2}(2\beta^2+H\beta_x)\Bigr]\\
 &=\frac{2\rho^2\beta^2+2H\rho^2\beta_x}{\tau^2\rho^2}
 =\frac2{\tau^2}\bigl(\beta^2+H\beta_x\bigr),
 \end{aligned}
\]
because $2(\tau^2-H^2)\beta^2+4H^2\beta^2=2\rho^2\beta^2$ and $2H\tau^2\beta_x+2H^3\beta_x=2H\rho^2\beta_x$. (Equivalently, $\Delta_0\psi=(H^2)_{xx}/\tau^2$.) By \eqref{eq:degenerate2},
\begin{equation}\label{eq:betaxx}
 \beta^2+H\beta_x\equiv0\quad\text{on }\mathcal D,\qquad\text{and then \eqref{eq:betay} gives}\quad \beta_y\equiv0 .
\end{equation}

\emph{Step 5: contradiction.} If $\beta(x_0)=0$ at some $x_0\in \mathcal D$, then $H_y(x_0)=-H\beta/\tau=0$, so $dH(x_0)=0$, contradicting $\mathcal D\subset\Omega$. Hence $\beta\ne0$ on $\mathcal D$. By \eqref{eq:betaxx}, $H=0$ would force $\beta^2=0$; so $H\ne0$ and $H_y=-H\beta/\tau\ne0$ on $\mathcal D$. Differentiating $\beta^2+H\beta_x=0$ in $y$ and using $\beta_y=\beta_{xy}=0$ gives $H_y\beta_x=0$, hence $\beta_x=0$ and then $\beta^2=0$, a contradiction.

Thus $\Omega$ is empty and $H$ is constant. The converse follows from \eqref{eq:holoidentity}.
\end{proof}

\section{Proof of Theorem~\ref{thm:main} for \texorpdfstring{$k\ne0$}{k not 0}}\label{sec:puiseuxsetup}

\begin{remark}[Why the round argument does not extend]\label{rem:whynot}
For $k=0$ the normalization $\QAR\equiv1$ fixes $p$ and $\lambda$ as functions
of $H$. For $k\ne0$ the term $-kA^2$, with $|A|^2=\lambda(1-u^2)/4$, brings in
$u$ as a second unknown, $\QAR$ may vanish on $\Omega$, and the normalization
determines neither $p$ nor $\lambda$; equivalently, $d\sigma=-ku\,dA\ne0$ in
\eqref{eq:dsigma}. We therefore work with the coordinates $(H,u)$.
\end{remark}

We first outline the argument for $k\ne0$. Suppose that
$\mathcal Q_{\mathrm{AR}}$ is holomorphic and $\Omega\ne\emptyset$. By
Lemma~\ref{lem:density} we may work near a point of $\Omega_*$, where $(H,u)$
are coordinates.
\begin{enumerate}
\item[(i)] In \S\ref{sec:resultant} the frame equations of
Lemma~\ref{lem:normal} become a first-order rational system
\eqref{eq:zH}--\eqref{eq:wu}, \eqref{eq:wH} for two functions $z,w$ of
$(H,u)$.
\item[(ii)] Equality of the mixed derivatives of $w$ yields a polynomial
relation $\Pi(H,u,w,z)=0$, and hence also $X\Pi=0$ (Lemma~\ref{lem:compat}).
\item[(iii)] In \S\ref{sec:algebraic}: every smooth solution of the system is
real analytic (Lemma~\ref{lem:analytic}; analyticity is a consequence, not a
hypothesis), so $z$ and $w$ generate, over $F:=\C(H,u)$, a field $E$ of germs.
\item[(iv)] Eliminating $z$ between $\Pi$ and $X\Pi$ by a resultant gives a
nontrivial polynomial equation for $w$ over $F$
(Lemmas~\ref{lem:resultant} and~\ref{lem:algebraic}); hence $E$ is a finite
extension of $F$.
\item[(v)] In \S\ref{sec:branches} we embed $E$ in a field of Puiseux series
in $u$, preserving $\partial_H$ and $\partial_u$ (Lemma~\ref{lem:extension}),
and show by comparing orders at $u=0$ that the system has no solution there
(Lemma~\ref{lem:puiseux}).
\item[(vi)] In \S\ref{sec:conclusion} these steps are combined to prove
Theorem~\ref{thm:main} for $k\ne0$.
\end{enumerate}

Fix $q\ne0$ and work on $\Omega_*$, where
$H\ne0$, $0<|u|<1$, $P\ne0$, and $(H,u)$ are coordinates. Put
\[
 s=1+H^2,\qquad \ell=4s+q,\qquad C=\frac1{1-u^2}+\frac q{4s}.
\]
On a simply connected neighborhood choose the angle $\theta$ and set
\begin{equation}\label{eq:zdef}
 \alpha:=\arctan H-\theta,\qquad z:=\tan\alpha,\qquad w:=s\,\alpha_H.
\end{equation}
With $\varphi=\arctan H$, so that $\cos\varphi=1/\sqrt s$ and
$\sin\varphi=H/\sqrt s$, the addition formulas and the definitions of $P,L$
(with $\tau=1$) give
\begin{equation}\label{eq:PLalpha}
 P=\sqrt s\cos\alpha,\qquad L=\sqrt s\sin\alpha,\qquad z=\tan\alpha=\frac LP,
\end{equation}
a real function, since $P\ne0$ on $\Omega_*$.

On $\Omega_*$ we have $H\ne0$, $0<|u|<1$, $q\ne0$, $1+z^2\ge1$, and $w\ne0$ by
\eqref{eq:wsratio} below; thus $(H,u,z,w)$ lies in $V:=\{Huw(1-u^2)\ne0\}$, and
all divisions below are legitimate.

\subsection{The rational system and a polynomial relation}\label{sec:resultant}

The angular equations give
\begin{align}
 z_H&=(1+z^2)\frac ws,\label{eq:zH}\\
 z_u&=(1+z^2)\mathcal A,\qquad \mathcal A:=-\frac{qHu}{4sw}+uzC,\label{eq:zu}\\
 w_u&=\mathcal G,\qquad \mathcal G:=\frac wu+\frac{qu(Hz+1)}{4s}-uwC(1+z^2).\label{eq:wu}
\end{align}
Indeed, $e_1=f\partial_H$ and $e_2=-2RP\,\partial_u$, as both sides agree on
$H$ and $u$ by \eqref{eq:eu}. By \eqref{eq:etheta} and \eqref{eq:ab} (with
$\tau=1$, $k=q$), $e_1\alpha=f/s+a=-quRP/(2s)$, so
$\alpha_H=-quRP/(2sf)$ and
\begin{equation}\label{eq:wsratio}
 w=s\alpha_H=-\frac{quRP}{2f},\qquad\text{i.e.}\qquad \frac f{RP}=-\frac{qu}{2w};
\end{equation}
then $z_H=(1+z^2)\alpha_H$ is \eqref{eq:zH}. Next,
$\alpha_u=e_2\alpha/(-2RP)=(-b+2uL/R)/(2RP)$; substituting $b$ from \eqref{eq:ab},
\eqref{eq:wsratio}, $L/P=z$ and $R^2=1-u^2$ gives $\alpha_u=\mathcal A$, hence
\eqref{eq:zu}. Finally, $w=-quR\sqrt s\cos\alpha/(2f)$ by \eqref{eq:PLalpha}, so
$w_u/w=1/u-u/(1-u^2)-z\alpha_u-f_u/f$; since $e_2f=af$, $f_u/f=-a/(2RP)
=-qu/(4sw)+qu/(4s)$, and substituting gives \eqref{eq:wu}.

To determine $w_H$, use $(\alpha_u)_H=(\alpha_H)_u$.
Introduce
\begin{equation}\label{eq:XYdef}
 X:=\partial_H+(1+z^2)\frac ws\,\partial_z+\mathcal B\,\partial_w,\qquad Y:=\partial_u+(1+z^2)\mathcal A\,\partial_z+\mathcal G\,\partial_w
\end{equation}
where $\mathcal B=w_H$. Along a solution, these vector fields act as
the total derivatives with respect to $H,u$.

Differentiating $\alpha_u=\mathcal A$ along the graph gives
\[
 (\alpha_u)_H=\mathcal A_H+\mathcal A_z\,z_H+\mathcal A_w\,w_H=\mathcal A_H+\mathcal A_z(1+z^2)\frac ws+\mathcal A_w\,\mathcal B.
\]
Since $(\alpha_H)_u=(w/s)_u=\mathcal G/s$, equality of mixed
derivatives gives
\begin{equation}\label{eq:wH}
\begin{split}
 w_H=\mathcal B&:=\frac{\mathcal G/s-\mathcal A_H-\mathcal A_z(1+z^2)w/s}{\mathcal A_w},\\
 \mathcal A_w&=\frac{qHu}{4sw^2}\ne0.
\end{split}
\end{equation}
Here subscripts on $\mathcal A$ are formal partial derivatives. With
$C_H=-qH/(2s^2)$, $\mathcal A_z=uC$ and
$\mathcal A_H=-qu(1-H^2)/(4s^2w)+uzC_H$, a direct simplification of
\eqref{eq:wH}, dividing term by term by $\mathcal A_w$, gives
\begin{equation}\label{eq:Bexplicit}
 \mathcal B=\frac{4w^3}{qHu^2}+\frac{w^2(3Hz+1)+w(1-H^2)}{Hs}-\frac{8w^3C(1+z^2)}{qH}.
\end{equation}

From now on, $X$ and $Y$ denote the vector fields \eqref{eq:XYdef} with $\mathcal B$ the explicit rational function \eqref{eq:Bexplicit}.

\begin{remark}[Geometric meaning of the system]\label{rem:dictionary}
With $\tau=1$, $\sigma=-z\,du/(1-u^2)$, and \eqref{eq:zH}--\eqref{eq:zu} are the
two components of $d\theta=u\sigma-\omega$ in \eqref{eq:dsigma}, i.e.\ of the
transport law \eqref{eq:T}; \eqref{eq:wu} is the symmetry $e_2f=af$ of the
Hessian of $H$; and \eqref{eq:wH} comes from $d(d\theta)=0$, i.e.\ the Gauss
equation. Thus $\Pi=0$ below expresses the compatibility of the Gauss equation
with the symmetry of the Hessian of $H$ for a surface whose normal is
$D$-parallel along $\nabla H$.
\end{remark}

We now derive a polynomial relation and use it to prove that
$w,z$ are algebraic over $\C(H,u)$.

Define the polynomials
\[
 \delta:=6s\ell-qu^2(14s+q)+q^2u^4
\]
and
\begin{align}
 \Pi:={}&4\delta w^3(1+z^2)\notag\\
 &+q\bigl[7u^2(\ell-qu^2)z^2+4s(11u^2-4)+7qu^2(1-u^2)\bigr]w^2\notag\\
 &-8Hq(1-u^2)w(3sw+qu^2)z\notag\\
 &-2q^2u^2(1-u^2)(1-3H^2)w-q^2su^2(1-u^2)\label{eq:Pi}
\end{align}
in $\C[H,u,w,z]$.

\begin{lemma}[Polynomial compatibility]\label{lem:compat}
Every solution of \eqref{eq:zH}--\eqref{eq:wH} satisfies $\Pi=X\Pi=0$.
\end{lemma}
\begin{proof}
On the graph $(H,u,z(H,u),w(H,u))$, $X,Y$ are the pushforwards of
$\partial_H,\partial_u$. Hence $[X,Y]$ vanishes there.

The $\partial_H$- and $\partial_u$-components of $[X,Y]$ vanish
identically. We calculate the other two components.

Put $\Xi=(1+z^2)w/s$. The $\partial_z$-component is
\[
 X\bigl((1+z^2)\mathcal A\bigr)-Y\Xi=(1+z^2)\Bigl(\mathcal A_H+\Xi\mathcal A_z+\mathcal B\mathcal A_w-\frac{\mathcal G}s\Bigr)+2z\Xi\mathcal A-2z\Xi\mathcal A.
\]
The terms $2z\Xi\mathcal A$ cancel, and the remaining bracket is zero
by \eqref{eq:wH}. Thus this component vanishes identically.

For the $\partial_w$-component, write
$\mathcal A=a_{-1}w^{-1}+a_0$, $\mathcal G=g_1w+g_0$,
and $\mathcal B=b_3w^3+b_2w^2+b_1w$, where
\[
\begin{aligned}
 a_{-1}&=-\frac{qHu}{4s},& a_0&=uzC,\\
 g_1&=\frac1u-uC(1+z^2),& g_0&=\frac{qu(Hz+1)}{4s},
\end{aligned}
\]
\[
 b_3=\frac4{qHu^2}-\frac{8C(1+z^2)}{qH},\qquad b_2=\frac{3Hz+1}{Hs},\qquad b_1=\frac{1-H^2}{Hs}.
\]
Put $v=(1+z^2)/s$ and $\tilde a_j=(1+z^2)a_j$ for $j=-1,0$.
No negative powers of $w$ occur in $X\mathcal G-Y\mathcal B$:
the only possible pole, in $\mathcal A$, is canceled by the factor
$w$ in $\partial_z\mathcal B$.

Expanding $X\mathcal G=\partial_H\mathcal G+vw\partial_z\mathcal G+\mathcal B\partial_w\mathcal G$
and $Y\mathcal B=\partial_u\mathcal B+(\tilde a_{-1}w^{-1}+\tilde a_0)\partial_z\mathcal B
+\mathcal G\partial_w\mathcal B$ in powers of $w$, one finds
$X\mathcal G-Y\mathcal B=\sum_{j=0}^3r_jw^j$, where
\begin{align}
 r_3&=-(b_3)_u-(1+z^2)a_0(b_3)_z-2g_1b_3,\notag\\
 r_2&=v(g_1)_z-(b_2)_u-(1+z^2)\bigl(a_{-1}(b_3)_z+a_0(b_2)_z\bigr)\notag\\
 &\quad{}-g_1b_2-3g_0b_3,\label{eq:rjformulas}\\
 r_1&=(g_1)_H+v(g_0)_z-(b_1)_u-(1+z^2)\bigl(a_{-1}(b_2)_z+a_0(b_1)_z\bigr)\notag\\
 &\quad{}-2g_0b_2,\notag\\
 r_0&=(g_0)_H-(1+z^2)a_{-1}(b_1)_z-g_0b_1.\notag
\end{align}

Substituting the explicit rational functions $a_{-1},a_0,g_1,g_0,b_3,b_2,b_1$ recorded above into the four formulas \eqref{eq:rjformulas} and simplifying---a direct but lengthy application of the product and quotient rules to functions already on hand, carried out in full in Fact~\ref{fact:rj}---gives
\[
\begin{gathered}
 r_0=-\frac{qu}{4Hs},\qquad r_1=-\frac{qu(1-3H^2+4Hz)}{2Hs^2},\\
 r_2=-\frac{J_2}{4Hus^2(u^2-1)},\qquad r_3=-\frac{\delta(1+z^2)}{qHus^2(u^2-1)},
\end{gathered}
\]
where
\[
 J_2:=7u^2(\ell-qu^2)z^2+4s(11u^2-4)+7qu^2(1-u^2)-24Hs(1-u^2)z.
\]
Multiplying by $-4qHus^2(u^2-1)$, the four terms $r_jw^j$ become
\[
\begin{gathered}
 4\delta(1+z^2)w^3,\qquad qJ_2w^2,\\
 -2q^2u^2(1-u^2)(1-3H^2)w-8Hq^2u^2(1-u^2)zw,\\
 -q^2su^2(1-u^2),
\end{gathered}
\]
respectively. The two terms containing $z$ to the first power, namely
$-24Hqs(1-u^2)zw^2$ from $qJ_2w^2$ and $-8Hq^2u^2(1-u^2)zw$, add up to
$-8Hq(1-u^2)w(3sw+qu^2)z$. Comparing with \eqref{eq:Pi} term by term,
the sum is exactly $\Pi$, so
\begin{equation}\label{eq:XG-YB}
 X\mathcal G-Y\mathcal B=-\frac{\Pi}{4qHus^2(u^2-1)}.
\end{equation}
On the graph, the left side is $w_{uH}-w_{Hu}=0$ and the
denominator is nonzero. Hence $\Pi=0$; differentiating this identity
along $X$ gives $X\Pi=0$.
\end{proof}

\subsection{Analyticity and algebraicity}\label{sec:algebraic}

We next pass from smooth solutions to analytic germs.

\begin{lemma}[Analyticity of the graph]\label{lem:analytic}
Every smooth graph solving
$z_H=(1+z^2)w/s$, $z_u=(1+z^2)\mathcal A$,
$w_H=\mathcal B$, $w_u=\mathcal G$ is locally real analytic on
$V=\{Huw(1-u^2)\ne0\}$.
\end{lemma}
\begin{proof}
The rational vector fields $X,Y$ are real analytic on $V$,
since $q\ne0$ and $s=1+H^2>0$.

Fix $p_0=(H_0,u_0,z_0,w_0)$ on the graph. The local flows
$\Phi_X^t,\Phi_Y^r$ are analytic in time and initial data
\cite[Ch.~I, \S1]{IlyashenkoYakovenko2008}, and, $X$ and $Y$ being locally
Lipschitz, every $C^1$ integral curve coincides with the corresponding flow
line. Hence
\[
 \Psi(h,v):=\Phi_Y^{v-u_0}\bigl(\Phi_X^{h-H_0}(p_0)\bigr)
\]
is analytic near $(H_0,u_0)$.

Since $(XH,Xu)=(1,0)$ and $(YH,Yu)=(0,1)$, the
$(H,u)$-projection of $\Psi(h,v)$ is $(h,v)$.

The graph curve at $u=u_0$ solves the $X$-equation through
$p_0$. Uniqueness gives
\[
 \Phi_X^{h-H_0}(p_0)=(h,u_0,z(h,u_0),w(h,u_0)),\qquad h\text{ near }H_0.
\]
For fixed $h$, the curve $v\mapsto(h,v,z(h,v),w(h,v))$ solves
the $Y$-equation with initial point $\Phi_X^{h-H_0}(p_0)$.
Uniqueness again gives $\Psi(h,v)=(h,v,z(h,v),w(h,v))$.

Thus $z,w$ are analytic. Uniqueness was applied only to integral curves lying on the
graph, so the relation $[X,Y]=0$ is not used.
\end{proof}

Fix a real point $(H_0,u_0)$ of the coordinate domain of $\Omega_*$. By
Lemma~\ref{lem:analytic}, $z$ and $w$ extend holomorphically to a complex
polydisc $\Delta$ around $(H_0,u_0)$ on which $H$, $u$, $w$, $1-u^2$ and $s$ have
no zeros; each identity in \eqref{eq:zH}--\eqref{eq:wH}, holding on the real
square $\Delta\cap\R^2$, holds on $\Delta$ by the identity theorem in each
variable. Inside the field of meromorphic germs at $(H_0,u_0)$ set
$E=\C(H,u)(w,z)$; $\C(H,u)$ embeds in it, since no nonzero polynomial vanishes
on an open set. All subsequent identities hold in this field.

We now eliminate $z$ by a resultant, to obtain a nontrivial relation
$\mathcal R(H,u,w)=0$.

Put $D_0:=Hqu^2s(u^2-1)$, $M:=D_0\mathcal B$ and $\mathcal N:=D_0\,X\Pi$, and
write $\Pi_0=\Pi|_{u=0}$, $\mathcal N_0=\mathcal N|_{u=0}$.

\begin{lemma}[Values at $u=0$]\label{lem:atzero}
The functions $M$ and $\mathcal N$ belong to $\C[H,u,w,z]$, with
\begin{equation}\label{eq:Npoly}
 \mathcal N=D_0\Pi_H+Hqu^2(u^2-1)(1+z^2)w\,\Pi_z+M\Pi_w,
\end{equation}
and
\begin{equation}\label{eq:Pi0N0}
 \Pi_0=8sw^2\Lambda(z),\qquad \Lambda(z)=3\ell w(1+z^2)-3Hqz-2q,\qquad
 \mathcal N_0=-4sw^3\,\partial_w\Pi_0 .
\end{equation}
\end{lemma}
\begin{proof}
Substituting $C=(4s+q(1-u^2))/(4s(1-u^2))$ in \eqref{eq:Bexplicit} and
multiplying by $D_0$ gives the polynomial
\[
\begin{aligned}
 M={}&4s(u^2-1)w^3+qu^2(u^2-1)\bigl[w^2(3Hz+1)+w(1-H^2)\bigr]\\
 &+2u^2\bigl[4s+q(1-u^2)\bigr]w^3(1+z^2).
\end{aligned}
\]
Since $X=\partial_H+(1+z^2)\frac ws\partial_z+\mathcal B\partial_w$ and
$D_0(1+z^2)w/s=Hqu^2(u^2-1)(1+z^2)w$, the formula for $\mathcal N$ follows,
and it is a polynomial. As $\deg_z\Pi=2$ and $\deg_zM=2$, it has
$\deg_z\mathcal N\le\max(2,\,2+1,\,2+2)=4$. At $u=0$: $\delta(0)=6s\ell$, so the five groups in
\eqref{eq:Pi} become $24s\ell w^3(1+z^2)$, $-16qsw^2$, $-24Hqsw^2z$, $0$ and
$0$, whose sum is $8sw^2\Lambda(z)$. In the formula for $\mathcal N$ the
coefficients $D_0$ and $Hqu^2(u^2-1)(1+z^2)w$ vanish at $u=0$, while
$M|_{u=0}=-4sw^3$; since setting $u=0$ commutes with $\partial_w$, this gives
$\mathcal N_0=-4sw^3\partial_w\Pi_0$.
\end{proof}

By \eqref{eq:Pi} and the proof of Lemma~\ref{lem:atzero}, $\deg_z\Pi=2$ and
$\deg_z\mathcal N\le4$. Let $\mathcal R:=\operatorname{Res}_z^{(2,4)}(\Pi,\mathcal N)$
be the Sylvester determinant of $\Pi$ and $\mathcal N$ formed with the nominal
degrees $(2,4)$ in $z$.

\begin{lemma}[Nonvanishing of the resultant]\label{lem:resultant}
The resultant $\mathcal R$ is a nonzero element of $\C[H,u,w]$, and
$\mathcal R=c_1\Pi+c_2\mathcal N$ for some $c_1,c_2\in\C[H,u,w,z]$.
\end{lemma}
\begin{proof}
We work in $\C[H,u,w,z]$, with $q\in\C^*$ fixed. The Sylvester determinant of
two polynomials of nominal degrees $(2,4)$ in $z$ lies in the ideal they
generate, whether or not the leading coefficients vanish
\cite[Ch.~IV, \S8]{Lang2002}, \cite[Ch.~3]{CLO2015}; this gives $c_1,c_2$. Setting
$u=0$ commutes with the determinant, so
$\mathcal R|_{u=0}=\operatorname{Res}_z^{(2,4)}(\Pi_0,\mathcal N_0)$, with
$\Pi_0,\mathcal N_0$ given by \eqref{eq:Pi0N0}.

Write $\Lambda(z)=3\ell w\,(z-\mu_1)(z-\mu_2)$ over $\overline{\C(H,w)}$; then
$\Pi_0$ has exact degree $2$ in $z$, with leading coefficient $a_2=24s\ell w^3$,
and for nominal degrees $(2,4)$ the root formula
$\operatorname{Res}_z^{(2,4)}(\Pi_0,\mathcal N_0)=a_2^{\,4}\,\mathcal N_0(\mu_1)\mathcal N_0(\mu_2)$
holds whatever the actual degree of $\mathcal N_0$ \cite[Ch.~IV, \S8]{Lang2002}.
Since $\partial_w\Pi_0=16sw\Lambda+24s\ell w^2(1+z^2)$, we get
$\mathcal N_0(\mu_j)=-96s^2\ell w^5(1+\mu_j^2)$, and
$(1+\mu_1^2)(1+\mu_2^2)=\Lambda(i)\Lambda(-i)/(3\ell w)^2
=q^2(4+9H^2)/(9\ell^2w^2)$, because $\Lambda(\pm i)=-q(2\pm3iH)$. Hence
\[
 \mathcal R|_{u=0}=2^{22}\,3^4\,q^2s^8\ell^4(4+9H^2)\,w^{20}\ne0,
\]
as $q\ne0$ and $s$, $\ell=4+4H^2+q$, $4+9H^2$ are nonzero polynomials in the
indeterminate $H$.
\end{proof}

\begin{lemma}[Algebraicity]\label{lem:algebraic}
For a real solution with $w\ne0$, $w$ is algebraic over $\C(H,u)$ and $z$ is
algebraic over $\C(H,u,w)$. Hence $[E:\C(H,u)]<\infty$.
\end{lemma}
\begin{proof}
By Lemma~\ref{lem:resultant}, $\mathcal R=c_1\Pi+c_2\mathcal N$ with polynomial
$c_1,c_2$. On the graph $\Pi=0$ and $\mathcal N=D_0\,X\Pi=0$, so
$\mathcal R(H,u,w(H,u))=0$ on an open polydisc. Consider $\deg_w\mathcal R$.
\emph{Case $\deg_w\mathcal R=0$.} Then $\mathcal R$ is a nonzero element of
$\C[H,u]$ vanishing on an open polydisc in the independent coordinates
$(H,u)$, which is impossible.
\emph{Case $\deg_w\mathcal R\ge1$.} Then $\mathcal R(H,u,Z)$ is a nonzero
polynomial in $Z$ over $\C(H,u)$, and $\mathcal R(H,u,w)=0$ shows that $w$ is
algebraic over $\C(H,u)$.

For $z$, collect powers of $z$ in $\Pi$: $\Pi=A_2z^2+A_1z+A_0$ with
\begin{align*}
 A_2&=4\delta w^3+7qu^2(\ell-qu^2)w^2,\\
 A_1&=-8Hq(1-u^2)w(3sw+qu^2),\\
 A_0&=4\delta w^3+q\bigl[4s(11u^2-4)+7qu^2(1-u^2)\bigr]w^2\\
 &\quad-2q^2u^2(1-u^2)(1-3H^2)w-q^2su^2(1-u^2).
\end{align*}
If $A_1\ne0$ we are done. If $A_1=0$, then $Hq(1-u^2)w\ne0$ gives
$w=-qu^2/(3s)$, an element of $F$; then $A_0$ is an element of $F$ as well,
and we show that it is not zero by expanding it in powers of $u$ in
$\C(H)((u))$. Since $\delta(0)=6s\ell$, the four terms of $A_0$ have, with
$w=-qu^2/(3s)$, lowest powers $u^6$, $u^4$, $u^4$ and $u^2$ respectively; the
last term is $-q^2su^2(1-u^2)$, with coefficient $-q^2s\ne0$ of $u^2$. Hence
$A_0=-q^2su^2+O(u^4)\ne0$.
Thus $\Pi$ is a nonzero polynomial in $z$ over $\C(H,u,w)$. If its
coefficients of positive degree in $z$ both vanished, then $\Pi=A_0\ne0$,
contradicting $\Pi=0$; hence $\Pi$ has positive degree in $z$, and $\Pi=0$
shows that $z$ is algebraic over $\C(H,u,w)$. Finiteness of $[E:\C(H,u)]$ follows, since $E$ is
generated by the two algebraic elements $w,z$.
\end{proof}

\subsection{Exclusion of the algebraic branches}\label{sec:branches}

Put $F=\C(H,u)$ and let $\mathbb K=\overline{\C(H)}$ be an algebraic
closure of $\C(H)$. The derivation $d/dH$ of
$\C(H)$ extends uniquely to $\mathbb K$ (Fact~\ref{fact:extension}(a)); we
denote it by a prime. Its field of constants is $\C$ (Fact~\ref{fact:extension}(b)):
\begin{equation}\label{eq:constants}
 \ker\bigl(d/dH|_{\mathbb K}\bigr)=\C.
\end{equation}
Let
\begin{equation}\label{eq:puiseuxfield}
 \mathcal L:=\bigcup_{n\ge1}\mathbb K\bigl((u^{1/n})\bigr)
\end{equation}
be the field of Puiseux series $\sum_{r\ge r_0}a_ru^r$, $a_r\in\mathbb K$,
$r\in\tfrac1n\Z$. Since $\mathbb K$ is algebraically closed of characteristic
zero, $\mathcal L$ is algebraically closed by the Newton--Puiseux theorem
\cite{Puiseux1850}, \cite[Ch.~IV]{Walker1950}. This is the reason for passing
from $\C(H)$ to $\mathbb K$: with coefficients in $\C(H)$ the corresponding
field is not algebraically closed (for instance, $Z^2=H$ has no root in
$\bigcup_{n}\C(H)((u^{1/n}))$), and the embedding of Lemma~\ref{lem:extension}
might not exist. For $0\ne\Phi\in\mathcal L$ let
$\nu(\Phi)$ be the smallest exponent occurring in $\Phi$, and $\nu(0)=+\infty$;
then $\nu(\Phi\Psi)=\nu(\Phi)+\nu(\Psi)$ and $\nu(\Phi+\Psi)\ge\min\{\nu(\Phi),\nu(\Psi)\}$,
with equality if $\nu(\Phi)\ne\nu(\Psi)$.

Define on $\mathcal L$
\begin{equation}\label{eq:DHDu}
 D_H\Bigl(\sum_ra_ru^r\Bigr):=\sum_r a_r'\,u^r,\qquad D_u\Bigl(\sum_ra_ru^r\Bigr):=\sum_r ra_r\,u^{r-1}.
\end{equation}
Both are derivations (the product rule holds coefficientwise). Under the
Laurent expansion at $u=0$, $F\hookrightarrow\C(H)((u))\subset\mathcal L$,
they restrict to $\partial_H,\partial_u$ on $F$: if $j:F\hookrightarrow\mathcal L$
denotes this embedding, then $D_H\circ j$ and $j\circ\partial_H$ (and likewise
$D_u\circ j$ and $j\circ\partial_u$) are derivations $F\to\mathcal L$ that vanish
on $\C$ and agree on $H$ and $u$, hence agree on $\C[H,u]$ and, by the quotient
rule, on $F$.
Directly from \eqref{eq:DHDu}: if $\nu(\Phi)=r$, write $\Phi=au^r(1+\eta)$
with $a\in\mathbb K^\times$ and $\nu(\eta)>0$; then $D_u\eta$ has order
$>-1$, and
\begin{equation}\label{eq:orders}
 \nu(D_H\Phi)\ge r,\qquad \frac{D_u\Phi}{\Phi}=\frac ru+\frac{D_u\eta}{1+\eta},\qquad \nu\Bigl(\frac{D_u\eta}{1+\eta}\Bigr)>-1 .
\end{equation}
In particular, $\nu(D_u\Phi)=r-1$ if $r\ne0$.

Since $F$ has characteristic zero, every algebraic extension of $F$ is
separable, and by Fact~\ref{fact:extension}(a) the derivations
$\partial_H,\partial_u$ of $F$ extend uniquely to any algebraic extension.
The next lemma states that an $F$-embedding intertwines these unique
extensions.

\begin{lemma}[Transport of the derivations]\label{lem:extension}
Let $E\supset F$ be a finite extension on which $\partial_H,\partial_u$ act as
derivations extending those of $F$, and let $\iota:E\hookrightarrow\mathcal L$ be
an $F$-embedding (it exists since $\mathcal L$ is algebraically closed). Then
\[
 D_H(\iota a)=\iota(\partial_Ha),\qquad D_u(\iota a)=\iota(\partial_ua)\qquad(a\in E).
\]
\end{lemma}
\begin{proof}
Let $\mathcal M\in F[Z]$ be the monic minimal polynomial of $a$; it is separable
(characteristic zero), so $\mathcal M_Z(a)\ne0$ and, $\iota$ being injective and
$F$-linear, $\mathcal M_Z(\iota a)=\iota(\mathcal M_Z(a))\ne0$. Let $\partial$ be
$\partial_H$ or $\partial_u$, with counterpart $D$ on $\mathcal L$, and let
$\mathcal M^\partial$ be $\mathcal M$ with $\partial$ applied to its coefficients.
Applying $\partial$ to $\mathcal M(a)=0$ in $E$, and $D$ to $\mathcal M(\iota a)=0$
in $\mathcal L$ (where $D=\partial$ on $F$), gives
\[
 \partial a=-\frac{\mathcal M^\partial(a)}{\mathcal M_Z(a)},\qquad D(\iota a)=-\frac{\mathcal M^\partial(\iota a)}{\mathcal M_Z(\iota a)}=\iota\Bigl(-\frac{\mathcal M^\partial(a)}{\mathcal M_Z(a)}\Bigr)=\iota(\partial a).\qedhere
\]
\end{proof}

Apply this to $E=\C(H,u)(w,z)$, which is finite over $F$ by
Lemma~\ref{lem:algebraic} and stable under $\partial_H,\partial_u$ by
\eqref{eq:zH}--\eqref{eq:wH}, whose denominators are products of powers of
$H$, $u$, $w$, $s$, $1-u^2$, all invertible in $E$ (germs not vanishing
identically on the connected $\Delta$ are invertible). The derivations of $E$
are the unique extensions of those of $F$ (Fact~\ref{fact:extension}(a)).
Identifying $E$ with $\iota(E)$, \eqref{eq:zH}--\eqref{eq:wH} and $\Pi=X\Pi=0$
hold in $\mathcal L$ with $D_H,D_u$, and $w\ne0$, $1+z^2\ne0$ persist.

\begin{remark}[The passage to $u=0$ is formal]\label{rem:formal}
From now on $u$ is an indeterminate of $F=\C(H,u)$, and ``$u=0$'' refers only to
the valuation $\nu$ of $\mathcal L$. No point of the surface with
$\langle N,\xi\rangle=0$ is used, the surface is not extended, and no
convergence of Puiseux series is claimed. Since $H,u$ are independent
indeterminates, $H$ is transcendental over $\C$; this makes the final identities
\eqref{eq:finalcase1} and \eqref{eq:finalcase2} contradictory.
\end{remark}

\begin{lemma}[Incompatibility of the orders]\label{lem:puiseux}
Let $q\in\C\setminus\{0\}$. There are no $w,z\in\mathcal L$ with $w\ne0$,
$1+z^2\ne0$, satisfying \eqref{eq:zH}--\eqref{eq:wu} (with $D_H,D_u$) and
$\Pi=0$.
\end{lemma}
\begin{proof}
All orders refer to the valuation $\nu$ of $\mathcal L$
(Remark~\ref{rem:formal}). Write $w=tu^\gamma+\cdots$ with $t\in\mathbb K^*$,
$\gamma=\nu(w)\in\Q$. The cases are summarized in Table~\ref{tab:puiseux}; in
each, the leading terms are compared exactly, including possible cancellations.

\begin{table}[ht]
\centering
\small
\renewcommand{\arraystretch}{1.3}
\begin{tabular}{|>{\raggedright\arraybackslash}p{1.9cm}|>{\raggedright\arraybackslash}p{2.8cm}|>{\raggedright\arraybackslash}p{2.8cm}|>{\raggedright\arraybackslash}p{3.2cm}|}
\hline
Case & Orders forced & Contradiction & Equation used; nonzero leading coefficient\\
\hline
$\nu(z)\ge0$, $\gamma=0$ & --- & term $t/u$ in \eqref{eq:wu} unmatched & \eqref{eq:wu}; $t\ne0$ ($t\in\mathbb K^*$)\\
$\nu(z)\ge0$, $\gamma\ge2$ & --- & $\nu(\Pi)=2$, so $\Pi\ne0$ & $\Pi$; coefficient $-q^2s$ of $u^2$ in the last term of \eqref{eq:Pi} ($q\ne0$, $s=1+H^2\ne0$)\\
$\nu(z)\ge0$, $0<\gamma<2$ & $\gamma=1$, $z=c+z_1u+\cdots$, $c\in\C$, $c\ne\pm i$ & order $2$ of $\Pi$ and order $1$ of \eqref{eq:zH} give \eqref{eq:finalcase1} & \eqref{eq:wu} ($\gamma=1$, as $t\ne0$), $\Pi$, \eqref{eq:zH}; $z_1\ne0$, $\varpi\ne0$ by \eqref{eq:order2}; constant term $3$ in \eqref{eq:finalcase1}\\
\hline
$\nu(z)=-b<0$, $\gamma>2$ & --- & term of order $1-\gamma<-1$ in \eqref{eq:logidentity} unmatched & \eqref{eq:logidentity}; coefficient $q/(4st)\ne0$ of $u^{1-\gamma}$\\
$\nu(z)=-b<0$, $\gamma<2$ & $\gamma-b=1$, $\gamma+b=2$, so $(\gamma,b)=(\tfrac32,\tfrac12)$ & $\nu(\Pi)=2$, so $\Pi\ne0$ & \eqref{eq:logidentity}, \eqref{eq:zu}, $\Pi$; coefficient $-q^2s$ of $u^2$ in the last term of \eqref{eq:Pi}\\
$\nu(z)=-b<0$, $\gamma=2$ & $b=2$, cancellation of the leading terms of $\mathcal A$ & orders $-1$ of \eqref{eq:logidentity} and $-2$ of \eqref{eq:zH} give \eqref{eq:finalcase2} & \eqref{eq:zu}: $D_uz$ has leading coefficient $-2c\ne0$, forcing the cancellation; $H^4$-coefficient $4$ in \eqref{eq:finalcase2}\\
\hline
\end{tabular}
\caption{The cases in the proof of Lemma~\ref{lem:puiseux}.}\label{tab:puiseux}
\end{table}

\emph{Preliminary: $\gamma\ge0$, and $\gamma\ge b$ if $\nu(z)=-b<0$.}
By \eqref{eq:zH}, $w=s\,D_Hz/(1+z^2)$. If $\nu(z)=-b<0$, then
$\nu(D_Hz)\ge-b$ and $\nu(1+z^2)=-2b$, so $\gamma\ge b$. If $\nu(z)\ge0$,
write $z=c+\varepsilon$ with $c\in\mathbb K$, $\nu(\varepsilon)>0$. When
$1+c^2\ne0$, $\nu(1+z^2)=0$ and $\gamma\ge\nu(D_Hz)\ge0$. When $c=\pm i$
(a constant), $D_Hz=D_H\varepsilon$ and $1+z^2=\varepsilon(2c+\varepsilon)$
has order exactly $\nu(\varepsilon)$, while $\nu(D_H\varepsilon)\ge\nu(\varepsilon)$;
again $\gamma\ge0$.

\medskip
\emph{Case 1: $\nu(z)\ge0$, $z=c+\cdots$ with $c\in\mathbb K$.}

\emph{Step 1: $\gamma=1$.} If $\gamma=0$, write $w=t+\eta$, $\nu(\eta)>0$. In
\eqref{eq:wu}, $D_uw=D_u\eta$ has order $>-1$, whereas $w/u$ has the term $t/u$
and the other two terms of $\mathcal G$ have order $\ge1$: a contradiction.
If $\gamma\ge2$, the first four terms of $\Pi$ (those containing $w$) have
orders at least $3\gamma,\ 2\gamma,\ \min\{2\gamma,\gamma+2\},\ \gamma+2$, all $>2$,
while $-q^2su^2(1-u^2)$ has order exactly $2$; so $\Pi\ne0$. Hence $0<\gamma<2$.
Then the three terms of $\mathcal G$ have orders $\gamma-1$, $\ge1$, $\ge\gamma+1$,
and comparing the terms of order $\gamma-1$ in \eqref{eq:wu} gives
$\gamma t=t$, that is, $\gamma=1$.

\emph{Step 2: expansions.} The right side of \eqref{eq:zH} has order $\ge\gamma=1$,
so the constant term of $D_Hz$ vanishes: $c'=0$, and $c\in\C$ by
\eqref{eq:constants}. Below, $o(1)$ denotes terms of positive order. Since
$w/u=t+o(1)$,
\[
 \mathcal A=-\frac{qH}{4s\,(w/u)}+uzC=-\frac{qH}{4st}+o(1).
\]
If $c=\pm i$, then $\varepsilon:=z-c\ne0$ has $\nu(D_u\varepsilon)=\nu(\varepsilon)-1$,
whereas $(1+z^2)\mathcal A=\varepsilon(2c+\varepsilon)\mathcal A$ has order $\nu(\varepsilon)$,
contradicting \eqref{eq:zu}. So $1+c^2\ne0$, and \eqref{eq:zu} gives
$D_uz=z_1+o(1)$ with
\[
 z_1=-\frac{(1+c^2)qH}{4st}\ne0 .
\]
As $D_u(z-c)$ has no terms of negative order, $z-c$ has no terms of order in $(0,1)$,
and its coefficient of $u$ is $z_1$:
$z=c+z_1u+(\text{terms of order}>1)$.

\emph{Step 3: coefficient identities.} Insert $w=tu+(\text{order}>1)$ and the expansion of $z$ into $\Pi$. The terms of order $2$ come from only three terms: $-16qst^2$
from the second group (the bracket equals $-16s+O(u^2)$), $-24Hqst^2c$ from the third
(only $3sw^2z$ matters), and $-q^2s$ from the last. Hence
\begin{equation}\label{eq:order2}
 -qs\bigl[8t^2(2+3Hc)+q\bigr]=0.
\end{equation}
Comparing the coefficients of $u$ in \eqref{eq:zH} gives $z_1'=(1+c^2)t/s$;
with the value of $z_1$ and $c'=0$,
\begin{equation}\label{eq:compare17}
 \Bigl(-\frac{qH}{4st}\Bigr)'=\frac ts .
\end{equation}
Set $\varpi=2+3Hc$. By \eqref{eq:order2}, $t^2=-q/(8\varpi)$ and $\varpi\ne0$;
differentiating $t^2\varpi=-q/8$ gives $t'/t=-3c/(2\varpi)$, hence
$(st)'/(st)=2H/s-3c/(2\varpi)$ and
\[
 \Bigl(-\frac{qH}{4st}\Bigr)'=-\frac q{4st}\Bigl[1-\frac{2H^2}s+\frac{3Hc}{2\varpi}\Bigr].
\]
Multiplying \eqref{eq:compare17} by $-4st/q$ and using $-4t^2/q=1/(2\varpi)$,
\[
 1-\frac{2H^2}s+\frac{3Hc}{2\varpi}=\frac1{2\varpi}.
\]
Multiplying by $2\varpi s$ and expanding,
\begin{equation}\label{eq:finalcase1}
 3-5H^2+(9H-3H^3)c=0 .
\end{equation}
Since $c\in\C$ and $H$ is transcendental over $\C$, this is a polynomial
identity in $H$ with constant term $3\ne0$, a contradiction. This excludes Case~1.

\medskip
\emph{Case 2: $\nu(z)=-b<0$, $z=cu^{-b}+\cdots$, $c\in\mathbb K^*$.}
By the preliminary, $\gamma\ge b>0$. Divide \eqref{eq:wu} by $w$ and multiply
\eqref{eq:zu} by $z/(1+z^2)$; adding, the terms $\pm qHuz/(4sw)$ cancel and
$-uC(1+z^2)+uz^2C=-uC$, so
\begin{equation}\label{eq:logidentity}
 \frac{D_uw}w+\frac{zD_uz}{1+z^2}=\frac1u-uC+\frac{qu}{4sw}.
\end{equation}
The left side is $\tfrac12D_u\Phi/\Phi$ with $\Phi=w^2(1+z^2)$, $\nu(\Phi)=2\gamma-2b$;
by \eqref{eq:orders} it equals $(\gamma-b)/u+(\text{order}>-1)$. On the right,
$-uC$ has order exactly $1$ (as $C(0)=\ell/(4s)\ne0$), and $qu/(4sw)$ has order $1-\gamma$.
If $\gamma>2$, the latter is a term of order $<-1$ that nothing on the left can
match. So $\gamma\le2$.

For $\gamma>b$, the first term of $\mathcal A$ has order $1-\gamma$ and the second
$1-b$, so $\nu(\mathcal A)=1-\gamma$, and comparing orders in \eqref{eq:zu}
($\nu(D_uz)=-b-1$) gives
\begin{equation}\label{eq:gammab}
 -b-1=1-\gamma-2b,\qquad\text{i.e.}\qquad \gamma+b=2 .
\end{equation}

\emph{Subcase $\gamma<2$.} Then $qu/(4sw)$ has order $>-1$, and the coefficients of
$u^{-1}$ in \eqref{eq:logidentity} give $\gamma-b=1$. So $\gamma>b$, and
\eqref{eq:gammab} yields $b=\tfrac12$, $\gamma=\tfrac32$. The five terms of $\Pi$
in \eqref{eq:Pi} then have orders $\tfrac72$, $3$, $\ge\tfrac52$, $\tfrac72$ and
exactly $2$ (the bracket in the second term has constant term $-16s$), so
$\nu(\Pi)=2$ and $\Pi\ne0$, a contradiction.

\emph{Subcase $\gamma=2$.} If $b<2$, then $\gamma>b$, so \eqref{eq:gammab} gives $b=0$,
which is excluded; since $b\le\gamma=2$, $b=2$. Write $w=tu^2+\cdots$,
$z=cu^{-2}+\cdots$. Both terms of $\mathcal A$ have order $-1$, so
$(1+z^2)\mathcal A$ has order $-5$ unless their leading coefficients cancel,
whereas $\nu(D_uz)=-3$. Cancellation requires
\[
 -\frac{qH}{4st}+\frac{\ell c}{4s}=0,\qquad\text{i.e.}\qquad c=\frac{qH}{\ell t}.
\]
In \eqref{eq:logidentity} the coefficient of $u^{-1}$ is $\gamma-b=0$ on the left and
$1+q/(4st)$ on the right, so $t=-q/(4s)$ and $c=-4sH/\ell$.
The coefficient of $u^{-2}$ in \eqref{eq:zH} gives $c'=c^2t/s=-4qH^2/\ell^2$.
On the other hand, with $s'=2H$ and $\ell'=8H$,
\[
 c'=-\frac{4\bigl[(1+3H^2)\ell-8H^2s\bigr]}{\ell^2}.
\]
Equating, $(1+3H^2)\ell-8H^2s=qH^2$; expanding with $\ell=4+4H^2+q$ this reads
\begin{equation}\label{eq:finalcase2}
 4s^2+q(1+2H^2)=0 .
\end{equation}
Its $H^4$-coefficient is $4\ne0$, impossible because $H$ is transcendental over $\C$.
All cases are excluded.

\end{proof}

\subsection{Conclusion of the proof}\label{sec:conclusion}

\begin{proof}[Proof of Theorem~\ref{thm:main} for $k\ne0$]
If $\Omega$ were nonempty, Lemma~\ref{lem:density} would give a point of
$\Omega_*$ with coordinates $(H,u)$, near which $z,w$ satisfy
\eqref{eq:zH}--\eqref{eq:wH} and $\Pi=X\Pi=0$, with $w\ne0$ by
\eqref{eq:wsratio} and $1+z^2=\sec^2\alpha\ge1$. Lemmas~\ref{lem:analytic},
\ref{lem:algebraic} and~\ref{lem:extension} transport these identities to
$\mathcal L$, contradicting Lemma~\ref{lem:puiseux}. Hence $dH=0$;
connectedness makes $H$ constant. The converse follows from
\eqref{eq:holoidentity}.
\end{proof}

\section{The product spaces}\label{sec:products}

In this section $\tau=0$, $k=\kappa$, and $\xi=\partial_t$ is parallel. We
write $\QAR$ for the same formal expression \eqref{eq:QAR} with $\tau=0$, that
is, $\QAR=2Hp-\kappa A^2$; the product spaces lie outside the scope of
Theorem~\ref{thm:main}. Data $(\lambda,u,H,p,A)$ on a simply connected domain
satisfying \eqref{eq:fund1}--\eqref{eq:fund2} define a local immersion into
$M^2(\kappa)\times\R$: for $\kappa\ne0$ these are the integrability
conditions (C.1)--(C.4) of \cite[Theorem~2.3, p.~283]{FM2007}, which are
sufficient when $\kappa-4\tau^2\ne0$ (see also \cite[Theorem~4.3, p.~101]{Daniel2007});
for $\kappa=0$ one adds the Gauss equation and the vertical equation
$A_\zeta-(\lambda_\zeta/\lambda)A=up$ (the pairing of
$\nabla_{\partial_\zeta}T$ with $\partial_\zeta$ in \eqref{eq:T}, using
\eqref{eq:christoffelconformal}) and uses the classical fundamental theorem
for surfaces in $\R^3$ \cite[\S4-3 and the Appendix to Ch.~4]{doCarmo1976};
the vertical equations make $T+uN$ a parallel unit field, identified with
$\partial_t$ by a rigid motion. Surfaces with constant angle function in
$S^2\times\R$ and $\mathrm H^2\times\R$ were classified in \cite{DFVV2007}
and \cite{DM2009}; we single out those with holomorphic $\calQAR$. Among
surfaces of constant angle $u_0\notin\{0,\pm1\}$, holomorphy of $\calQAR$ is
equivalent to $\nabla H$ being parallel to $T$ on $\{dH\ne0\}$, that is, to $H$
being constant along the lines of curvature orthogonal to $T$
(Corollary~\ref{cor:constantangle}). This is the condition, additional to the
data of the classifications in \cite{DFVV2007,DM2009}, that selects the family
of Proposition~\ref{prop:constantangle}: there all data depend only on the
parameter $\sigma$ along $T$.

Fix $\kappa\in\R$, $u_0\in(-1,1)\setminus\{0\}$, and put $R_0=\sqrt{1-u_0^2}$.
Let $Y(\sigma)>0$ be a solution of $(\log Y)''=-\kappa u_0^2Y^2$ on an interval;
this equation has the first integral
\begin{equation}\label{eq:Yode}
 \Bigl(\frac{Y'}{Y}\Bigr)^2=C-\kappa u_0^2Y^2 ,\qquad C\in\R\ \text{constant}.
\end{equation}
On the strip $\zeta=\sigma+i\upsilon$ put $v=\log Y$ and
\begin{equation}\label{eq:constantangledata}
 \lambda=Y^2,\qquad u=u_0,\qquad A=\frac{R_0Y}2,\qquad H=\frac{R_0v'}{2u_0Y},\qquad p=-\frac{\lambda H}2 .
\end{equation}

\begin{proposition}[Constant angle surfaces]\label{prop:constantangle}
The data \eqref{eq:constantangledata} define a local immersion into
$M^2(\kappa)\times\R$ with constant angle function $u_0$,
\begin{equation}\label{eq:HY}
 \QAR\equiv-\frac{R_0^2C}{4u_0^2},\qquad H'=-\frac{R_0C}{2u_0Y},\qquad
 4H^2+\kappa R_0^2=\frac{R_0^2C}{u_0^2Y^2}.
\end{equation}
In particular, $\calQAR$ is holomorphic, and the immersion is CMC if and only if
$C=0$. Moreover:
\begin{enumerate}
\item[(i)] if $\kappa>0$, then $C>0$;
\item[(ii)] if $\kappa=0$, then $C\ge0$, and $C=0$ exactly when $Y$ is constant
and $H=p=0$ (a plane of constant angle);
\item[(iii)] if $\kappa<0$, every sign of $C$ occurs; the surfaces with $C<0$
satisfy $4H^2<|\kappa|R_0^2$.
\end{enumerate}
\end{proposition}
\begin{proof}
By hypothesis $v''=-\kappa u_0^2Y^2$, and \eqref{eq:Yode} holds with a
constant $C$ because $(v'^2+\kappa u_0^2e^{2v})'=2v'(v''+\kappa u_0^2e^{2v})=0$.
All data depend on $\sigma$ only, so $F_\zeta=F_{\bar\zeta}=F'/2$, with
$'=d/d\sigma$, and we check \eqref{eq:fund1}--\eqref{eq:fund2}, the Gauss
equation and the vertical equation.
\begin{itemize}
\item $A_{\bar\zeta}=\tfrac{u\lambda}2H$, i.e.\ $A'=u_0\lambda H$: both sides equal $R_0Yv'/2$.
\item $p_{\bar\zeta}=\tfrac\lambda2(H_\zeta+\kappa uA)$, i.e.\ $p'=\lambda(H'/2+\kappa u_0A)$:
$H'=\frac{R_0}{2u_0Y}(v''-v'^2)=-\frac{R_0C}{2u_0Y}$, and
$p'=-\frac{R_0}{4u_0}(Yv')'=-\frac{R_0Y}{4u_0}(v'^2+v'')
=-\frac{R_0YC}{4u_0}+\frac{\kappa u_0R_0Y^3}2=\lambda\bigl(\tfrac{H'}2+\kappa u_0A\bigr)$.
\item $u_\zeta=-HA-\frac{2p}\lambda\bar A$: both sides vanish, since $2p/\lambda=-H$ and $A$ is real.
\item $4|A|^2=R_0^2Y^2=\lambda(1-u_0^2)$.
\item Gauss: $\det S=H^2-4|p|^2/\lambda^2=0$, and $(\log\lambda)''=2v''=-2\kappa u_0^2\lambda$,
which is Liouville's formula for $K=\kappa u_0^2$.
\item Vertical equation: $A_\zeta-(\lambda_\zeta/\lambda)A=\frac{R_0Yv'}4-\frac{R_0Yv'}2=-\frac{R_0Yv'}4=u_0p$.
\end{itemize}
Finally
$\QAR=2Hp-\kappa A^2=-\lambda H^2-\kappa R_0^2Y^2/4
=-\frac{R_0^2}{4u_0^2}(v'^2+\kappa u_0^2Y^2)=-\frac{R_0^2C}{4u_0^2}$, and
the last identity in \eqref{eq:HY} follows from $4H^2=R_0^2v'^2/(u_0^2Y^2)$ and \eqref{eq:Yode}.
If $\kappa>0$, \eqref{eq:Yode} gives $C\ge\kappa u_0^2Y^2>0$. If $\kappa=0$,
then $C=v'^2\ge0$, and $C=0$ if and only if $Y$ is constant, in which case
$H=p=0$. If $\kappa<0$ and $C>0$, examples are given by
Proposition~\ref{prop:products} below; for $C<0$, by homothety we may take
$\kappa=-1$, and $Y=\sqrt{|C|}/(|u_0|\cos(\sqrt{|C|}\,\sigma))$ on
$|\sqrt{|C|}\,\sigma|<\pi/2$ satisfies $Y'/Y=\sqrt{|C|}\tan(\sqrt{|C|}\sigma)$ and
$(Y'/Y)^2=|C|/\cos^2(\sqrt{|C|}\sigma)-|C|=u_0^2Y^2+C$; finally $C<0$ in
\eqref{eq:HY} gives $4H^2<|\kappa|R_0^2$.
\end{proof}

Fix $\mathfrak a\ne0$, put $\mathfrak m=\sqrt{1+\mathfrak a^2}$, and let
$(y(\sigma),r(\sigma))$, with $y>0$, solve
\begin{equation}\label{eq:ydef}
 y'=-\mathfrak a y\,r,\qquad r'=\mathfrak a\kappa y^2,\qquad r^2+\kappa y^2=1 .
\end{equation}
The constraint is preserved, since $(r^2+\kappa y^2)'=2\mathfrak a\kappa y^2r-2\mathfrak a\kappa y^2r=0$.
For $\kappa\le0$, $r^2\ge1$, so $r$ has a fixed sign, which may be either; for
$\kappa>0$, $r$ may change sign.
On the strip $\zeta=\sigma+i\upsilon$ consider the data
\begin{equation}\label{eq:productdata}
 u=\frac{\mathfrak a}{\mathfrak m},\qquad \lambda=\mathfrak m^2y^2,\qquad H=-\frac{r}{2\mathfrak m y},\qquad
 p=-\frac{\lambda H}2=\frac{\mathfrak m yr}4,\qquad A=\frac y2 .
\end{equation}

\begin{proposition}[The hypothesis $\tau\ne0$ is essential]\label{prop:products}
For every $\kappa\in\R$, the data \eqref{eq:productdata} define a local
immersion into $M^2(\kappa)\times\R$ with $\QAR\equiv-1/4$ and
$H'=-\mathfrak a/(2\mathfrak my)\ne0$. In particular, $\calQAR$ is holomorphic
and the mean curvature is not constant.
\end{proposition}
\begin{proof}
These are the surfaces of Proposition~\ref{prop:constantangle} with
$u_0=\mathfrak a/\mathfrak m$, $R_0=1/\mathfrak m$, $Y=\mathfrak my$ and
$C=\mathfrak a^2$. Indeed $v'=y'/y=-\mathfrak ar$, so
$(Y'/Y)^2=\mathfrak a^2(1-\kappa y^2)=C-\kappa u_0^2Y^2$; then
$A=R_0Y/2=y/2$, $H=R_0v'/(2u_0Y)=-r/(2\mathfrak my)$,
$\QAR=-R_0^2C/(4u_0^2)=-1/4$ and $H'=-R_0C/(2u_0Y)=-\mathfrak a/(2\mathfrak my)\ne0$.
\end{proof}

Explicit solutions of \eqref{eq:ydef}, with $x:=\mathfrak a\sigma+\sigma_0$, are:
$y=y_0e^{-\mathfrak a\sigma}$, $r=1$ ($y_0>0$, all $\sigma$) for $\kappa=0$;
$y=\kappa^{-1/2}\operatorname{sech}x$, $r=\tanh x$, for all $x\in\R$, when $\kappa>0$;
and $y=|\kappa|^{-1/2}\operatorname{csch}x$, $r=\coth x$, on $x>0$, when $\kappa<0$.
For $\kappa>0$ the mean curvature $H=-\sqrt\kappa\sinh x/(2\mathfrak m)$ vanishes
exactly on the curve $x=0$, where $r=0$, and $H'\ne0$ there.
For $\kappa=\pm1$, restricted to $x>0$, this is, up to notation, the family of
\cite[Example~3.1, p.~286]{FM2007}, where the original differential
$P\,dz^2=-\tfrac14dz^2$ is observed to be holomorphic although $H$ is not
constant; for other $\kappa\ne0$ it follows by a homothety. The examples are
therefore known; what is new is that constant angle is necessary and that the
family of Proposition~\ref{prop:constantangle} exhausts all possibilities
(Theorem~\ref{thm:products}). They have zero extrinsic curvature, as required by
\cite[p.~173]{CL2007}. Their normalized differential
$Q_{\mathrm{FM}}=\QAR/H=-1/(4H)$ is defined on $\{H\ne0\}$ and is not holomorphic
there; when $\kappa>0$ it does not extend continuously across the curve $H=0$,
since $\QAR=-1/4\ne0$, while $\calQAR$ is holomorphic on the whole surface. Thus,
at $\tau=0$, the original and normalized differentials have different
exceptional surfaces (for $Q_{\mathrm{FM}}$ in $\mathrm H^2\times\R$, the
rotational surfaces of \cite[Example~3.2]{FM2007}); the nonvanishing of $\tau$
is used in Lemma~\ref{lem:density} and in the round proof. By \eqref{eq:HY}, the
surfaces of Proposition~\ref{prop:products} are, after a change of conformal
parameter, exactly those of Proposition~\ref{prop:constantangle} with $C>0$ and
the same $u_0$ (end of Step~5 below); for $\kappa<0$, those with $C<0$ satisfy
$4H^2<|\kappa|R_0^2$ and are not among them.

In the following theorem, $\nabla H$, $u$ and $T$ are as in
\S\ref{sec:formalism}, now with $\tau=0$ and $k=\kappa$.

\begin{theorem}[Classification in the product spaces]\label{thm:products}
Let $\kappa\ne0$, and let $\Sigma$ be a surface in $M^2(\kappa)\times\R$ with
holomorphic $\calQAR$ and nonconstant mean curvature. Then $dH$ vanishes
nowhere on $\Sigma$, and:
\begin{enumerate}
\item[(a)] the angle function is a constant $u_0\notin\{0,\pm1\}$;
\item[(b)] $\nabla H$ is a principal direction with principal curvature $0$,
and $T$ is parallel to $\nabla H$;
\item[(c)] $|\nabla H|=\dfrac{|u_0|\,\bigl|4H^2+\kappa(1-u_0^2)\bigr|}{2\sqrt{1-u_0^2}}$;
\item[(d)] $\Sigma$ is locally congruent to a surface of
Proposition~\ref{prop:constantangle} with the same $u_0$ and with $C\ne0$ of the
sign of $4H^2+\kappa(1-u_0^2)$; up to local congruence it is determined by
$\kappa$, $u_0$ and this sign. When $\kappa>0$ or $4H^2>|\kappa|(1-u_0^2)$, it is, after a change of
conformal parameter, locally congruent to a surface of
Proposition~\ref{prop:products} with the same $u_0$, including near points where
$H=0$ if $\kappa>0$.
\end{enumerate}
\end{theorem}

\begin{proof}
\emph{Outline.} Steps~1--3 show that $du=0$ on a dense open subset of
$\Omega=\{dH\ne0\}$, Step~4 deduces (a)--(c), Step~5 identifies each component of
$\Omega$ with a model, and Step~6, the only global step, shows $\Omega=\Sigma$ by
Gr\"onwall's inequality.

\emph{Reduction.} Holomorphy is again equivalent to \eqref{eq:holo} with
$\tau=0$, and the proof of Proposition~\ref{prop:normalform} applies verbatim:
on $\Omega$, in the frame $e_1=\nabla H/f$, $e_2=Je_1$, $f=|\nabla H|$,
\[
 [S]_{(e_1,e_2)}=\begin{pmatrix}0&0\\0&2H\end{pmatrix},\qquad \det S=0,\qquad K=\kappa u^2 .
\]
This proves the first half of (b). The proof of Lemma~\ref{lem:normal} uses
$\tau\ne0$ only through $\rho^2=H^2+\tau^2\ne0$; on the open set
$\Omega_0:=\Omega\cap\{H\ne0,\ 0<|u|<1\}$ it therefore gives, with $\tau=0$,
$P=H\sin\theta$, $L=H\cos\theta$, $R=\sqrt{1-u^2}$ and
$\nabla_{e_1}e_1=a\,e_2$, $\nabla_{e_2}e_1=b\,e_2$ (here $a,b$ are the
connection coefficients of Lemma~\ref{lem:normal}),
\begin{equation}\label{eq:prodframe}
\begin{gathered}
 a=-\frac{\kappa uR\sin\theta}{2H},\qquad b=-\frac fH-\frac{\kappa uR\cos\theta}{2H},\\
 e_1u=0,\qquad e_2u=-2HR\sin\theta,\\
 e_1\theta=-a,\qquad e_2\theta=-b+\frac{2Hu\cos\theta}{R},\qquad e_2f=af .
\end{gathered}
\end{equation}
These formulas only require $H\ne0$ and $|u|<1$; in particular they hold
where $u=0$.

\emph{Step 1: the branch $\sin\theta\equiv0$.} Let $V$ be an open set with
$H\ne0$, $|u|<1$ and $\sin\theta\equiv0$. Then $\cos\theta=\varepsilon\in\{\pm1\}$
is locally constant, $T=\varepsilon Re_1$, $a=0$, and $e_1u=e_2u=0$, so $u$ is
locally constant on $V$. The equation for $e_2\theta$ in \eqref{eq:prodframe}
gives $b=2\varepsilon Hu/R$, and comparing with the expression of $b$,
\begin{equation}\label{eq:prodf}
 f=-\frac{\varepsilon u\,(4H^2+\kappa R^2)}{2R}.
\end{equation}
In particular $u\ne0$ on $V$, since $f>0$. The relation $e_2f=af=0$ holds
because $f$ is a function of $H$, and the Gauss equation is automatic:
$K=-e_1b-b^2=u^2(4H^2+\kappa R^2)/R^2-4H^2u^2/R^2=\kappa u^2$.

\emph{Step 2: $\Omega_0$ is dense in $\Omega$.} The sets $\{H=0\}$,
$\{|u|=1\}$ and $\{u=0\}$ are closed in $\Omega$; we show they have empty
interior. If $H\equiv0$ on an open set, then $dH=0$ there. If $|u|\equiv1$ on
an open set, then $T\equiv0$ and \eqref{eq:T} gives $uS\equiv0$, so
$S\equiv0$ and $H\equiv0$. If $u\equiv0$ on an open set, shrink it so that
$H\ne0$; there $R=1$ and $0=e_2u=-2H\sin\theta$, so $\sin\theta\equiv0$, and
Step~1 gives $u\ne0$, a contradiction.

\emph{Step 3: $\sin\theta=0$ on $\Omega_0$.} Suppose $\sin\theta\ne0$ on a
nonempty open set $U\subset\Omega_0$. There
$dH\wedge du(e_1,e_2)=-2fHR\sin\theta\ne0$, so $(H,u)$ are local coordinates,
with $e_1=f\partial_H$ and $e_2=-2HR\sin\theta\,\partial_u$. Put
$c=\cos\theta$ and $\mathfrak n=4H^2+\kappa R^2$. Dividing \eqref{eq:prodframe} by
$e_1H=f$ and $e_2u=-2HR\sin\theta$,
\begin{equation}\label{eq:prodsystem}
 f_u=\frac{\kappa u}{4H^2}\,f,\qquad c_H=-\frac{\kappa uR(1-c^2)}{2Hf},\qquad c_u=\frac{2fR+uc\mathfrak n}{4H^2R^2}.
\end{equation}
(For the first, $a/(-2HR\sin\theta)=\kappa u/(4H^2)$; for the last two,
$c_H=-\sin\theta\,\theta_H$ and $c_u=-\sin\theta\,\theta_u$.)

\emph{The Gauss equation gives $f_H$.} With $\omega=a\,e_1^\flat+b\,e_2^\flat$,
$d\omega=-K\,dA$ reads $K=e_2a-e_1b-a^2-b^2$. Using \eqref{eq:prodsystem},
\[
\begin{aligned}
 e_2a&=\kappa(1-2u^2)(1-c^2)-\frac{\kappa uc\,(2fR+uc\mathfrak n)}{4H^2},\\
 a^2&=\frac{\kappa^2u^2R^2(1-c^2)}{4H^2},\qquad b^2=\frac{(2f+\kappa uRc)^2}{4H^2},\\
 e_1b&=-\frac{ff_H}H+\frac{f^2}{H^2}+\frac{\kappa^2u^2R^2(1-c^2)}{4H^2}+\frac{\kappa uRcf}{2H^2}.
\end{aligned}
\]
Substituting in $K=\kappa u^2$ and collecting terms, where the terms in
$c^2\mathfrak n$ give $\kappa u^2c^2+\kappa^2u^2R^2c^2/(4H^2)$, we obtain
\begin{equation}\label{eq:prodfH}
 f_H=\frac{\mathcal E}{Hf},\qquad
 \mathcal E:=2f^2+2\kappa uRcf+\frac{\kappa^2u^2R^2}2+\kappa H^2\bigl(3u^2-1+R^2c^2\bigr).
\end{equation}

\emph{Compatibility of $f$.} From \eqref{eq:prodsystem},
$f_{uH}=\kappa uf_H/(4H^2)-\kappa uf/(2H^3)$, while \eqref{eq:prodfH} gives
$f_{Hu}=\mathcal E_u/(Hf)-\kappa u\mathcal E/(4H^3f)$. Equating and
multiplying by $4H^3f$,
\begin{equation}\label{eq:prodcompat}
 2H^2\mathcal E_u-\kappa u\mathcal E+\kappa uf^2=0 .
\end{equation}
Differentiating $\mathcal E$ with \eqref{eq:prodsystem}, $(uR)_u=(1-2u^2)/R$
and $(R^2)_u=-2u$, the terms in $c^2$ cancel and one finds
\[
 \mathcal E_u=\frac{2\kappa uf^2}{H^2}+3\kappa Rcf+\frac{\kappa^2u^2Rcf}{H^2}
 +\kappa^2u(1-2u^2)+6\kappa H^2u+\frac{\kappa^2uR^2c^2}2 .
\]
Substituting in \eqref{eq:prodcompat}, the terms in $c^2$ and in $u^2Rcf$
cancel, and after division by $\kappa/2\ne0$,
\begin{equation}\label{eq:prodPi}
 12H^2Rf\,c=-u\,\mathcal U,\qquad \mathcal U:=6f^2+24H^4+2\kappa H^2(3-7u^2)-\kappa^2u^2R^2 .
\end{equation}

\emph{Elimination of $c$.} Since $12H^2Rf\ne0$, \eqref{eq:prodPi} gives $c$.
Differentiating it in $u$, with $f_u$ from \eqref{eq:prodsystem},
\[
 \Bigl(\frac u{Rf}\Bigr)_u=\frac{4H^2-\kappa u^2R^2}{4H^2R^3f},\qquad
 \mathcal U_u=\frac{3\kappa uf^2}{H^2}-28\kappa H^2u-2\kappa^2u(1-2u^2),
\]
and comparing with the value of $c_u$ in \eqref{eq:prodsystem}, in which
\[
 uc\mathfrak n=-\frac{u^2\mathcal U\mathfrak n}{12H^2Rf},
\]
we obtain, after multiplication by $48H^4R^3f$ and division by $R^2$,
\[
 \mathcal U\,(2\kappa u^2-4H^2)-4H^2u\,\mathcal U_u-24H^2f^2=0 .
\]
Expanding, the terms $\pm12\kappa u^2f^2$ cancel, and
\begin{equation}\label{eq:prodf2}
\begin{aligned}
 24H^2f^2=\Phi:={}&-48H^6-12\kappa H^4+108\kappa H^4u^2+12\kappa^2H^2u^2\\
 &-24\kappa^2H^2u^4+\kappa^3u^4(u^2-1).
\end{aligned}
\end{equation}

\emph{Contradiction.} By \eqref{eq:prodsystem}, $(f^2)_u=\kappa uf^2/(2H^2)$,
so \eqref{eq:prodf2} gives $2H^2\Phi_u=\kappa u\Phi$ on $U$. For each fixed
$H$ this is an identity between polynomials in $u$ on an interval; the
right side has degree $7$ in $u$, with leading coefficient $\kappa^4$,
while the left side has degree at most $5$. As $\kappa\ne0$, this is
impossible. Hence $\sin\theta\equiv0$ on $\Omega_0$.

\emph{Step 4: conclusion of (a)--(c).} By Steps~1 and~3, $du=0$ on the dense
open set $\Omega_0$, hence on $\Omega$, and $u\equiv u_0$ on each component
of $\Omega$; by Step~2, $u_0\notin\{0,\pm1\}$. This is (a); $T=\varepsilon
R_0e_1$ gives the second half of (b), and \eqref{eq:prodf} with $f>0$ gives
(c) on $\Omega_0$, hence on $\Omega$ by continuity.

\emph{Step 5: congruence.} Fix a component $\Omega_1$ of $\Omega$, with
constant angle $u_0$. By (c), $\mathfrak n_0:=4H^2+\kappa R_0^2$ does not
vanish on $\Omega_1$, so it has a constant sign $\varsigma\in\{\pm1\}$, and
\eqref{eq:prodf} gives $\varepsilon=-\varsigma\operatorname{sign}u_0$. On
$\Omega_1$, $e_2H=0$ and $[e_1,e_2]=-b\,e_2$, where $b=2\varepsilon Hu_0/R_0$
and $f$ are functions of $H$ alone. Hence, near any point, there is a
coordinate $t$ with $e_1=f(H)\partial_H$ and $e_2=\mu(H)\partial_t$,
$\mu'/\mu=-b/f$, and the metric $dH^2/f^2+dt^2/\mu^2$, the shape operator
($Se_1=0$, $Se_2=2He_2$), $T=\varepsilon R_0e_1$ and $u=u_0$ are given, as
functions of $(H,t)$, by formulas depending only on $u_0$, $\varsigma$ and
$\kappa$ (the constant in $\mu$ is absorbed into $t$).

Now take a surface of Proposition~\ref{prop:constantangle} with the same
$u_0$ and with $\operatorname{sign}C=\varsigma$; by Proposition~\ref{prop:constantangle}
and \eqref{eq:HY} this is possible exactly when $\varsigma=1$ or $\kappa<0$, which
covers all cases allowed by (c). It has holomorphic $\calQAR$, constant angle
$u_0$, and $H'=-R_0C/(2u_0Y)\ne0$, so $H$ is strictly monotone in $\sigma$,
and Steps~1--4 apply to it with the same $u_0$ and $\varsigma$, hence the same
$\varepsilon$. Its mean curvature fills the whole interval of values allowed
by $\varsigma$. Indeed, on a maximal solution of $v''=-\kappa u_0^2e^{2v}$:
if $\kappa>0$ (so $C>0$), $Y$ increases to its maximum, where $v'=0$ and
$H=0$, and then decreases, with $Y\to0$ and $v'\to\mp\sqrt C$ at the ends, so
$H$ runs over all of $\R$; if $\kappa<0$ and $C>0$, $v'$ never vanishes,
$Y$ runs over $(0,\infty)$, and $4H^2=R_0^2(C/(u_0^2Y^2)+|\kappa|)$ runs over
$(|\kappa|R_0^2,\infty)$, with the sign of $H$ fixed by that of $v'/u_0$, which
can be chosen; if $\kappa<0$ and $C<0$, $Y$ has a minimum where $H=0$ and
tends to $\infty$ at both ends, so $4H^2$ runs over $[0,|\kappa|R_0^2)$ and $H$
over $(-h_0,h_0)$, $h_0=\sqrt{|\kappa|}R_0/2$. Hence the model surface contains
points with any value of $H$ taken on $\Omega_1$, and, in the coordinates
$(H,t)$ above, both surfaces have the same metric, shape operator, $T$ and $u$.
By the uniqueness part of the fundamental theorem \cite[Theorem~4.3, p.~101]{Daniel2007},
$\Omega_1$ is locally congruent to the model. When $\kappa>0$, or
$4H^2>|\kappa|R_0^2$, we have $C>0$; the model is then a surface of
Proposition~\ref{prop:products} with the same $u_0$. Indeed, put
$\mathfrak a:=u_0/R_0$, so that $\mathfrak m=1/R_0$ and $u_0=\mathfrak a/\mathfrak m$,
and change the conformal parameter by $\zeta=b\tilde\zeta$, $b:=|\mathfrak a|/\sqrt C$;
this is a change of coordinates on the surface, not a homothety of the ambient
space. Then $\tilde Y(\tilde\sigma):=bY(b\tilde\sigma)$ satisfies
$(\log\tilde Y)''=-\kappa u_0^2\tilde Y^2$ with constant $\tilde C=b^2C=\mathfrak a^2$
in \eqref{eq:Yode}. Setting $y:=\tilde Y/\mathfrak m$ and $r:=-\tilde v'/\mathfrak a$, with
$\tilde v=\log\tilde Y$, one checks $r^2=1-\kappa y^2$, $y'=-\mathfrak ayr$ and
$r'=\kappa u_0^2\tilde Y^2/\mathfrak a=\mathfrak a\kappa y^2$, so $(y,r)$ solves
\eqref{eq:ydef}, and the data \eqref{eq:constantangledata} become
\eqref{eq:productdata}, with $H=-r/(2\mathfrak my)$. No change of orientation is
needed, since $r$ is allowed to take either sign. This includes the points where
$H=0$, which occur exactly when $\kappa>0$.

Steps~1--5 establish (a)--(d) on each connected component of $\Omega$.

\emph{Step 6: $\Omega=\Sigma$.} Since $H$ is not constant, $\Omega\ne\emptyset$.
Suppose $\Omega\ne\Sigma$, and let $U$ be a component of $\Omega$, with
constant angle $u_0$ and $R_0=\sqrt{1-u_0^2}$. By (c),
\[
 |\nabla H|=c_0\,|n_0(H)|\ \text{ on }U,\qquad c_0=\frac{2|u_0|}{R_0},\quad n_0(H)=H^2+\frac{\kappa R_0^2}4 .
\]
Since $\Sigma$ is connected, there is a piecewise smooth path
$\gamma:[0,1]\to\Sigma$, with $|\gamma'|\le m$ for a constant $m>0$, such that
$\gamma(0)\in U$ and $\gamma(1)\notin U$. Let $s_*=\inf\{s:\gamma(s)\notin U\}$; as
$U$ is open, $\gamma([0,s_*))\subset U$ and $x:=\gamma(s_*)\notin U$. Then
$dH(x)=0$: otherwise $x\in\Omega$, and a connected neighborhood of $x$ in
$\Omega$ meets $U$, so $x\in U$. Letting $s\to s_*$, the continuity of $H$ and
$|\nabla H|$ gives $n_0(H(x))=0$. If $\kappa>0$ this is impossible. If
$\kappa<0$, then $H(x)=\epsilon_1h_0$ with $h_0=\sqrt{|\kappa|}R_0/2$ and
$\epsilon_1=\pm1$. Put $\chi(s)=H(\gamma(s))-\epsilon_1h_0$. The function $\chi$ is
Lipschitz, hence absolutely continuous, and differentiable except at the
finitely many break points of $\gamma$. Near $s_*$ we have
$|H(\gamma(s))+\epsilon_1h_0|\le3h_0$, hence
\[
\begin{aligned}
 |\chi'(s)|&\le m\,|\nabla H|(\gamma(s))=mc_0|\chi(s)|\,|H(\gamma(s))+\epsilon_1h_0|\\
 &\le M_0|\chi(s)|,\qquad M_0:=3mc_0h_0,
\end{aligned}
\]
on $[s_*-\eta_0,s_*)$ for some $\eta_0>0$, except at finitely many points. As
$\chi(s_*)=0$, absolute continuity gives
$|\chi(s)|\le\int_s^{s_*}M_0|\chi(r)|\,dr$, and Gr\"onwall's inequality gives
$\chi\equiv0$ on $[s_*-\eta_0,s_*]$. Thus $H=\epsilon_1h_0$ at the points
$\gamma(s)\in U$, $s\in[s_*-\eta_0,s_*)$, where (c) gives $|\nabla H|=0$,
contradicting $U\subset\Omega$. Hence $\Omega=\Sigma$, that is, $dH$ vanishes
nowhere, and $u\equiv u_0$ on $\Sigma$. No completeness is used.
\end{proof}

Theorem~\ref{thm:productsintro} is Theorem~\ref{thm:products}: item (d) of the
former follows from Step~5, where the metric, the shape operator, $T$ and $u$
are expressed in terms of $\kappa$, $u_0$ and the sign of $4H^2+\kappa R_0^2$ only.

\begin{proof}[Proof of Corollary~\ref{cor:closedproducts}]
If $\Sigma$ is closed, $H$ attains its maximum at a point where $dH=0$. By
Theorem~\ref{thm:products}, $H$ is constant. This agrees with
\cite[pp.~172--173]{CL2007}.
\end{proof}

\begin{corollary}\label{cor:constantangle}
Let $\kappa\ne0$ and let $\Sigma$ be a surface in
$M^2(\kappa)\times\R$ with constant angle function $u\equiv u_0\notin\{0,\pm1\}$.
Then $\calQAR$ is holomorphic if and only if $\nabla H$ is parallel to $T$ on
$\{dH\ne0\}$, that is, $H$ is constant along the lines of curvature
orthogonal to $T$.
\end{corollary}
\begin{proof}
Since $u$ is constant, \eqref{eq:vertical} with $\tau=0$ gives
$\langle SX,T\rangle=-X(u)=0$ for all $X$, so $ST=0$; as $T\ne0$, $T$ is a
principal direction with principal curvature $0$, and the other principal
direction $JT$ has principal curvature $2H$. On $\{dH\ne0\}$, the argument of
Proposition~\ref{prop:Dparallel} with $\tau=0$ (where $D=\overline\nabla$)
shows that holomorphy is equivalent, pointwise, to $S\nabla H=0$. Where
$H\ne0$, the kernel of $S$ is spanned by $T$, so $S\nabla H=0$ means
$\nabla H\parallel T$; since $\{H\ne0\}$ is dense in $\{dH\ne0\}$ and both
conditions are closed, holomorphy on $\{dH\ne0\}$ is equivalent to
$\nabla H\parallel T$ there. On the interior of $\{dH=0\}$, $H$ is locally
constant and $(\QAR)_{\bar\zeta}=0$ by \eqref{eq:holoidentity}; by continuity
$(\QAR)_{\bar\zeta}$ then vanishes everywhere. The converse is
part~(b) of Theorem~\ref{thm:products}.
\end{proof}

\begin{remark}\label{rem:H2R}
In the surfaces of Proposition~\ref{prop:constantangle} with $\kappa<0$, $H$
approaches $\pm h_0$ only where $Y\to\infty$, that is, at infinite distance,
in accordance with Step~6 of the proof of Theorem~\ref{thm:products}. For $\kappa=0$, $\calQAR=2Hp$, and holomorphy of
$Hp$ for surfaces in $\R^3$ is a classical problem outside the scope of this
paper.
\end{remark}

\section{Proof of Corollary~\ref{cor:zeros}}\label{sec:global}

By Theorem~\ref{thm:main}, $H$ is constant; compactness is used only to count
zeros. With isothermal parameters \cite{Chern1955}, $\Sigma$ is a compact
Riemann surface of genus $\mathsf g$, and since $p$ and $A$ transform as
coefficients of a $(2,0)$- and a $(1,0)$-form, $\calQAR$ is a holomorphic
section of $K_\Sigma^{\otimes2}$, of degree $4\mathsf g-4$ \cite[Ch.~III]{FK1980}.
If $\calQAR\not\equiv0$ it has $4\mathsf g-4$ zeros counted with multiplicity;
for $\mathsf g=0$ this is impossible, so $\calQAR\equiv0$, and the rotational
invariance follows from the references given after Corollary~\ref{cor:zeros}.
For $k=0$ it is direct: $p\equiv0$, so $\Sigma$ is a totally umbilical, hence
round, sphere of $S^3(\tau^2)$.
\qed

\appendix
\section{Auxiliary computations}\label{sec:preliminaries}

This appendix contains two facts used in Section~\ref{sec:puiseuxsetup}:
Fact~\ref{fact:rj} is the direct (but lengthy) computation of the coefficients in the proof of
Lemma~\ref{lem:compat}, and Fact~\ref{fact:extension} collects standard facts on derivations of
algebraic extensions.
\begin{fact}[Coefficients of the compatibility identity]\label{fact:rj}
With $a_{-1}$, $a_0$, $g_1$, $g_0$, $b_3$, $b_2$, $b_1$ as in the proof of Lemma~\ref{lem:compat}, the coefficients $r_0,\dots,r_3$ of \eqref{eq:rjformulas} are
\[
\begin{gathered}
 r_0=-\frac{qu}{4Hs},\qquad r_1=-\frac{qu(1-3H^2+4Hz)}{2Hs^2},\\
 r_2=-\frac{J_2}{4Hus^2(u^2-1)},\qquad r_3=-\frac{\delta(1+z^2)}{qHus^2(u^2-1)},
\end{gathered}
\]
with $J_2$ as defined after \eqref{eq:rjformulas} in the main text.
\end{fact}
\begin{proof}
Put $v:=(1+z^2)/s$, $V_0:=1+z^2$ and $d:=1-u^2$. Since $(b_1)_z=0$,
$r_0=(g_0)_H-g_0b_1$. The quotient rule gives
\begin{align*}
 (g_0)_H&=\frac{qu}4\cdot\frac{z\,s-(Hz+1)\cdot2H}{s^2}=\frac{qu}{4s^2}\bigl[z(s-2H^2)-2H\bigr]\\
 &=\frac{qu}{4s^2}\bigl[z(1-H^2)-2H\bigr].
\end{align*}
Since $g_0b_1=qu(Hz+1)(1-H^2)/(4Hs^2)$,
\begin{align*}
 r_0&=\frac{qu}{4Hs^2}\Bigl[H\bigl(z(1-H^2)-2H\bigr)-(Hz+1)(1-H^2)\Bigr]\\
 &=\frac{qu}{4Hs^2}\Bigl[Hz(1-H^2)-2H^2-Hz(1-H^2)-(1-H^2)\Bigr]\\
 &=\frac{qu}{4Hs^2}\bigl({-}s\bigr)=-\frac{qu}{4Hs},
\end{align*}
where $2H^2+(1-H^2)=s$.
For the other coefficients, $C_H=-qH/(2s^2)$, $C_u=2u/d^2$, and
$(b_1)_u=(b_1)_z=(b_2)_u=0$. First,
\begin{align*}
r_1&=(g_1)_H+v(g_0)_z-V_0a_{-1}(b_2)_z-2g_0b_2\\
&=\frac{quH V_0}{2s^2}+\frac{quH V_0}{4s^2}
  +\frac{3quH V_0}{4s^2}
  -\frac{qu(Hz+1)(3Hz+1)}{2Hs^2}\\
&=\frac{qu}{2Hs^2}\bigl[3H^2(1+z^2)-(Hz+1)(3Hz+1)\bigr]\\
&=-\frac{qu(1-3H^2+4Hz)}{2Hs^2}.
\end{align*}
Next, using $(b_3)_z=-16Cz/(qH)$ and $(b_2)_z=3/s$, the terms containing
$zCV_0/s$ are $v(g_1)_z=-2uzCV_0/s$, $-V_0a_0(b_2)_z=-3uzCV_0/s$ and
$-V_0a_{-1}(b_3)_z=-4uzCV_0/s$, whose sum is $-9uzCV_0/s$; the remaining terms
are $-g_1b_2$ and $-3g_0b_3$, expanded with $g_1=1/u-uCV_0$. Hence
\begin{align*}
r_2&=v(g_1)_z-V_0\bigl[a_{-1}(b_3)_z+a_0(b_2)_z\bigr]
      -g_1b_2-3g_0b_3\\
&=-\frac{9uzCV_0}{s}-\frac{3Hz+1}{Hsu}
  +\frac{uCV_0(3Hz+1)}{Hs}\\
&\hspace{1em}-\frac{3(Hz+1)}{Hsu}
  +\frac{6uCV_0(Hz+1)}{Hs}\\
&=\frac{7uCV_0-(6Hz+4)/u}{Hs}\\
&=\frac{7u^2(4s+qd)(1+z^2)-4sd(6Hz+4)}{4Hus^2d}.
\end{align*}
Here $4s+qd=\ell-qu^2$, and the constant term in the numerator is
$7u^2(4s+qd)-16sd=4s(11u^2-4)+7qu^2d$.
Finally, $(b_3)_u=-8/(qHu^3)-8C_uV_0/(qH)$ gives
\begin{align*}
r_3&=-(b_3)_u-V_0a_0(b_3)_z-2g_1b_3\\
&=\frac1{qH}\left[\frac8{u^3}+8C_uV_0+16uz^2C^2V_0\right.\\
&\hspace{3em}\left.-2\left(\frac1u-uCV_0\right)
                  \left(\frac4{u^2}-8CV_0\right)\right]\\
&=\frac1{qH}\left[8C_uV_0+\frac{24CV_0}u+16uC^2V_0\bigl(z^2-V_0\bigr)\right]\\
&=\frac{8V_0}{qH}\left(C_u+\frac{3C}{u}-2uC^2\right),
\end{align*}
since the terms $\pm8/(qHu^3)$ cancel, $-2\cdot\frac1u\cdot(-8CV_0)-2(-uCV_0)\frac4{u^2}=\frac{24CV_0}u$, and $z^2-V_0=-1$. Expanding the last bracket,
\begin{align*}
C_u+\frac{3C}{u}-2uC^2
&=\frac{2u}{d^2}+\frac3{ud}+\frac{3q}{4su}
 -\frac{2u}{d^2}-\frac{qu}{sd}-\frac{q^2u}{8s^2}\\
&=\frac{24s^2+6qs-14qsu^2-q^2u^2+q^2u^4}{8s^2ud}
 =\frac{\delta}{8s^2ud}.
\end{align*}
Collecting the four coefficients gives the stated formulas for $r_0,\dots,r_3$.
\end{proof}

\begin{fact}\label{fact:extension}
Let $F$ be a field of characteristic zero with a derivation $D$.
\begin{enumerate}
\item[(a)] If $a$ is algebraic over $F$ with minimal polynomial $\mathcal M$, then
$\mathcal M_Z(a)\ne0$ and $D$ extends uniquely to a derivation of $F(a)$, with
\[
 D(a)=-\frac{\mathcal M^D(a)}{\mathcal M_Z(a)},
\]
where $\mathcal M^D$ applies $D$ to the coefficients of $\mathcal M$. Consequently $D$
extends uniquely to any algebraic extension of $F$.
\item[(b)] For $F=\C(H)$ and $D=d/dH$, the field of constants of the extension to
$\mathbb K=\overline{\C(H)}$ is $\C$.
\end{enumerate}
\end{fact}
\begin{proof}
(a) This is standard \cite[Ch.~VIII, \S5]{Lang2002}. Applying $D$ to $\mathcal M(a)=0$ gives
$\mathcal M^D(a)+\mathcal M_Z(a)D(a)=0$, which forces the formula. Conversely, $\mathcal M_Z$
is invertible modulo $\mathcal M$ ($\mathcal M$ is irreducible and $\deg\mathcal M_Z<\deg\mathcal M$),
so the formula defines a derivation of $F[Z]$ into $F[Z]/(\mathcal M)$ that kills $\mathcal M$,
hence descends to $F(a)\cong F[Z]/(\mathcal M)$. Uniqueness at each finite stage gives
compatibility along towers, hence the extension to algebraic extensions.

(b) If $a/b\in\C(H)$ is reduced and $(a/b)'=0$, then $a'b=ab'$ and $a\mid a'$; as $\deg a'<\deg a$
unless $a$ is constant, $a$ is constant, and then $ab'=0$ makes $b$ constant. So the constants
of $\C(H)$ are $\C$. Now let $c\in\mathbb K$ with $c'=0$, with minimal polynomial
$\mathcal M=Z^m+\sum_{j<m}a_jZ^j$ over $\C(H)$. By (a), $0=\mathcal M^D(c)=\sum_{j<m}a_j'c^j$,
a polynomial relation of degree $<m$; by minimality all $a_j'=0$, so $a_j\in\C$. An
irreducible polynomial over $\C(H)$ with coefficients in $\C$ has degree $1$ ($\C$ being
algebraically closed), so $m=1$ and $c\in\C$.
\end{proof}

\section*{Data availability statement}

This manuscript has no associated data.

\section*{Conflict of interest statement}

The authors declare that they have no conflict of interest.

\end{document}